\documentclass[oneside,leqno,draft]{amsart}
\usepackage{latexsym,amssymb,amsmath,graphics}

\begin{document}
%
%

\def\ct{\mathcal{C}_{\rho}}
\def\bl1{\mathcal{B}_{\rho,\epsilon,\lambda}}
\def\el{{\epsilon, \lambda}}
\def\rel{{\rho,\epsilon,\lambda}}
\def\lel{{2\lambda,\epsilon,\lambda}}
\def\Sd{\mathbb{S}^{d-1}}
\newcommand{\Dis}{\mbox{Dis}}
\newcommand{\Lip}{\mbox{Lip}}
\renewcommand{\k}{\kappa}
\newcommand{\Dt}{\Delta_t}
\def\F{{\mathcal{F}}}
\newcommand{\supp}{{\rm supp}}
\def\bP{\overline{\Phi}}

\def\aa{{\mathcal{A}}}
\def\gg{{\mathcal{G}}}
\def\kk{{\mathcal{K}}}
\def\bb{{\mathcal{B}}}
\def\Ss{{\mathcal{S}}}
\def\nn{{\mathcal{N}}}
\def\A{{\mathbb{A}}}
\def\B{{\mathbb{B}}}
\def\C{{\mathbb{C}}}
\def\E{{\mathbb{E}}}
\def\H{{\mathbb{H}}}
\def\M{{\mathbb{M}}}
\def\m{{\mathbb{m}}}
\def\N{{\mathbb{N}}}
\def\P{{\mathbb{P}}}
\def\Q{{\mathbb{Q}}}
\def\R{{\mathbb{R}}}
\def\S{{\mathbb{S}}}
\def\T{{\mathbb{T}}}
\def\X{{\mathbb{X}}}
\def\Y{{\mathbb{Y}}}
\def\Z{{\mathbb{Z}}}
\def\PP{{\mathbb{P}}}

\def\w{{\mathrm{w}}}

\newcommand{\ol}{{\mathcal L}}
\newcommand{\on}{{\mathcal N}}

\newcommand{\jd}{\frac{1}{2}}
\newcommand{\pf}{\psi f}
\newcommand{\8}{\infty}
\renewcommand{\d}{\delta}
\renewcommand{\a}{\alpha}
\renewcommand{\b}{\beta}
\newcommand{\D}{\Delta}
\renewcommand{\O}{\Omega}
\renewcommand{\L}{\Lambda}
\newcommand{\p}[1]{\P\big[#1\big]}
\newcommand{\e}[1]{\E\big[#1\big]}
\newcommand{\eps}{\varepsilon}
\newcommand{\z}{\zeta}
\newcommand{\ba}{{\bf \underline a}}
\renewcommand{\bb}{{\bf \underline b}}
\newcommand{\bc}{{\bf \underline c}}
\newcommand{\g}{\gamma}
\newcommand{\s}{\sigma}
\newcommand{\cab}{c_{\ba,\bb}}
\newcommand{\ov}{\overline}
\newcommand{\wt}{\widetilde}
\newcommand{\cj}{\overline c_j}
\newcommand{\cp}{\overline c_p}
\newcommand{\pab}{P_{j}^{\ba,\bb}}
\newcommand{\nc}{\nn^\C}
\newcommand{\1}{{\bf 1}}
\renewcommand{\wt}{\widetilde}
\newcommand{\wh}{\widehat}
\newcommand{\ro}{\rho}
\renewcommand{\t}{\theta}
\renewcommand{\th}{\theta}
\renewcommand{\d}{\delta}
\renewcommand{\e}{\varepsilon}
\renewcommand{\o }{\omega}
\newcommand{\is}[2]{\langle #1,#2\rangle}
\newcommand{\iss}[1]{\langle #1\rangle}
\newcommand{\bis}[2]{\Big\langle #1,#2\Big\rangle}
\newcommand{\bbis}[2]{\bigg\langle #1,#2\bigg\rangle}
\newcommand{\suf}[1]{\lceil #1\rceil}
\newcommand{\ttv}{T_{t,v}}
\newcommand{\psitv}{\psi_{t,v}}
\newcommand{\etv}{\eta_{t,v}}
\newcommand{\dmn}{d\mu(h_1)\ldots d\mu(h_n)}
\newcommand{\chitx}{\chi_{\d_t(x)}}
\newcommand{\ep}{{\frac \e 2}}
\newcommand{\mar}[1]{\marginpar{\bf \footnotesize #1}}
\newcommand{\xm}{{x_-}}
\newcommand{\xp}{{x_+}}
\newcommand{\ym}{{y_-}}
\newcommand{\yp}{{y_+}}
\newcommand{\vm}{{V_{\frac \a2,-}}}
\newcommand{\vp}{{V_{\frac \a2,+}}}
\newcommand{\ttau}{{\widetilde \tau}}
\newcommand{\dt}{b_n}
\newcommand{\rp}{\R ^d\setminus \{0\}}
\newcommand{\ups}{\Upsilon}
\renewcommand{\l}{\lambda}
\newcommand{\nlc}{\nn_\l^\C}
\newcommand{\tel}{\th,\e,\l}

\newtheorem{problem}{Problem}

\newtheorem{con}[equation]{Condition}
\newtheorem{thm}[equation]{Theorem}
\newtheorem{cor}[equation]{Corollary}
\newtheorem{lem}[equation]{Lemma}
\newtheorem{lemma}[equation]{Lemma}
\newtheorem{mthm}[equation]{Main Theorem}
\newtheorem{mth}[equation]{Main Theorem}
\newtheorem{prop}[equation]{Proposition}
\theoremstyle{definition}
\newtheorem{defn}[equation]{Definition}
\newtheorem{rem}[equation]{Remark}
\newtheorem{ex}{Example}[section]

\numberwithin{equation}{section}


\title[Asymptotics of the stationary solution ]{Asymptotics of
stationary solutions of multivariate stochastic recursions with heavy
tailed inputs and related limit theorems}
\author[D. Buraczewski, E. Damek, M.Mirek]
{Dariusz Buraczewski, Ewa Damek, Mariusz Mirek}
\address{D. Buraczewski, E. Damek, M. Mirek \\
Uniwersytet Wroclawski\\
Instytut Matematyczny\\
pl. Grunwaldzki 2/4\\
 50-384 Wroclaw\\
 Poland}
\email{dbura@math.uni.wroc.pl\\ edamek@math.uni.wroc.pl\\ mirek@math.uni.wroc.pl}

\thanks{
D.~Buraczewski and E.~Damek were partially supported by MNiSW N N201 393937. M.~Mirek was partially
supported by MNiSW grant N N201 392337. D.~Buraczewski was also supported by European Commission via
 IEF
Project (contract number
PIEF-GA-2009-252318 - SCHREC)}

\begin{abstract}
Let $\Phi_n$ be an i.i.d. sequence of Lipschitz mappings of $\R^d$. We study the Markov  chain
$\{X_n^x\}_{n=0}^\8$ on $\R^d$ defined by the recursion $X_n^x = \Phi_n(X^x_{n-1})$, $n\in\N$,
$X_0^x=x\in\R^d$. We assume that $\Phi_n(x)=\Phi(A_n x,B_n(x))$ for a fixed continuous function
$\Phi:\R^d\times \R^d\to\R^d$, commuting with dilations and i.i.d random pairs $(A_n,B_n)$, where
$A_n\in {\rm End}(\R^d)$ and $B_n$ is a continuous mapping of $\R^d$. Moreover, $B_n$ is
$\a$-regularly varying and $A_n$ has a faster decay at infinity than $B_n$. We prove that the
stationary measure $\nu$ of the Markov chain $\{X_n^x\}$ is $\a$-regularly varying. Using this
result we show that, if $\a<2$, the partial sums $S_n^x=\sum_{k=1}^n X_k^x$, appropriately
normalized, converge to an $\a$-stable random variable. In particular, we obtain new results
concerning the random coefficient
autoregressive process $X_n = A_n X_{n-1}+B_n$.
\end{abstract}

\maketitle


\section{Introduction and main results}

We consider the vector space $\R^d$ endowed with an arbitrary norm
$|\cdot|$. We fix once for all  a continuous mapping
$\Phi:\R^d\times\R^d \to \R^d$, commuting with dilations, i.e.
$\Phi(tx,ty) = t\Phi(x,y)$ for every $t>0$. Let $(A, B)$ be a
random pair, where $A\in{\rm End}(\R^d)$ and $B$ is a continuous
mapping of $\R^d$. We assume that $B$ is of the form
$B(x)=B^1+B^2(x)$, where $B^1$ is a random vector in $\R^d$ and
$B^2$ is a random mapping of $\R^d$ such that $|B^2(x)|\le
B^3|x|^{\d_0}$ for every $x\in\R^d$, where $\d_0\in [0,1)$ is a
fixed number and  $B^3\ge 0$ is random.
Given a sequence $(A_n, B_n)_{n\in\N}$ of independent random copies of the generic pair $(A, B)$ and a starting point $x\in \R^d$,
we
define the Markov chain by
\begin{equation}
\label{gen_rec}
\begin{split}
X_0^x &= x,\\
X_n^x &  = \Phi(A_n X^x_{n-1}, B_n(X^x_{n-1})),\ \mbox{for $n\in\N$.}
\end{split}
\end{equation}
If $x=0$ we just write for simplicity $X_n$ instead of $X_n^{0}$.
Also, to simplify the notation, let $\Phi_n (x) = \Phi(A_n x,
B_n(x))$. Then the definition above can be expressed in a more
concise way, $X_n^x = \Phi_n(X_{n-1}^x)$.

\medskip

The main example we have in mind 
is a random coefficient
autoregressive process on $\R^d$, called also a random
difference equation or an affine stochastic recursion. This
process is defined by
\begin{equation}
\label{affine rec}
X_{1,n}^x = A_n X_{1,n-1}^x + B_n.
\end{equation}
and as one can easily see it is a particular example of
\eqref{gen_rec}, just by taking  $\Phi(x,y)=x+y$ and $B_n^2\equiv
0$.

For an another example take $d=1$, $\Phi (x,y)=\max (x,y)$ and $B_n^2 \equiv 0$. Then we obtain the random extremal equation
\begin{equation}\label{max rec}
 X_{2,n}^x = \max (A_n X_{2,n-1}^x, B_n),
\end{equation} studied e.g. by Goldie \cite{G}.

\medskip

In this paper we assume that the Markov chain $\{X_n^x\}$ is $\gamma$-geometric. This means
that there are constants $0<C<\8$ and $0<\rho<1$ such that the  moment
of order $\gamma>0$ of the Lipschitz coefficient of
$\Phi_{n}\circ\ldots\circ\Phi_{1}$ decreases exponentially fast as
$n$ goes to infinity, i.e.
\begin{equation}
\label{eq: geom1}
\E\Big[ \big| X_n^x - X_n^y \big|^\g \Big] \le C \rho^n|x-y|^{\g} ,\quad n\in\N,x,y\in\R^d.
\end{equation}
We say that a random vector $W\in \R^d$ is regularly varying with index $\a>0$
(or $\a$-regularly varying) if there is a slowly varying function
$L$ such that the limit
\begin{equation}
\label{regular_varying}
\lim_{t \to \infty}t^{\a }L(t)\E\big[ f(t^{-1} W)\big] =\int_{\rp}f(x)\L(dx)=: \langle f,\Lambda\rangle,
\end{equation}
exists for every $f\in C_c (\rp )$ and thus defines a Radon measure
$\L$ on $\rp$. The measure $\L$ will be called the tail measure.
It can be easily checked that
$\int_{\rp}f(rx)\L(dx) = r^{\a }\langle f,\L \rangle$ for every $r>0$, and so the tail measure $\L $ is $\a$-homogeneous, i.e. in radial coordinates we have
\begin{equation}
\label{radial}
\langle f,\L \rangle =\int_0^{\8}\int_{\mathbb{S}^{d-1}}f(r\o )\ \s_\L (d\o
)\frac{dr}{r^{1+\a }},
\end{equation}
for some measure $\s_\L$ on the unit sphere $\Sd \subseteq\R^d$. The measure $\s_\L$ will we called the spherical measure of $\L$. Observe that $\s_\L$ is nonzero if and only if $\L$
is nonzero.

\medskip

Under  mild assumptions there exists a unique stationary distribution $\nu$ of $\{X_n^x\}$ (see Lemma \ref{lem: geom}).
The main purpose of this paper is to prove, under some further hypotheses, that the distribution $\nu$
is $\a$-regularly varying and next to obtain a limit theorem for partial sums $S_n^x = \sum_{k=1}^nX_{k}^x$.

\medskip

Our first main result is the following
\begin{thm}
\label{thm: ogon} Let $\{X_n^x\}$ be the Markov chain defined by
\eqref{gen_rec}. Assume that
\begin{itemize}
\item $B^1$ is $\a$-regularly varying
with the nonzero tail measure $\L_b$ and the corresponding slowly varying function $L_b$
is bounded away from zero and infinity on any compact set;
\item the Markov chain $\{X_n^x\}$ is $\g$-geometric for some $\g>\a$;
\item there exists $\b>\a$ such that $\E\|A\|^\b <\8$;
\item there exists $\eps_0>0$ such that $\E\big[ (B^3)^{\frac{\a}{\d_0}+\eps_0} \big]<\8$, if $0<\d_0<1$ and  $\E\big[ (B^3)^{\a+\eps_0} \big]<\8$, if $\d_0=0$;
\item $\P[B^1: \Phi(0,B^1)\not=0]>0$.
\end{itemize}
Then  the Markov chain $\{X_n^x\}$ has a unique stationary measure $\nu$. If $X$ is a random variable distributed according to $\nu$, then
 $X$  is $\a$-regularly varying with a nonzero tail measure $\L^1$, i.e. for every $f\in C_c (\rp )$
\begin{equation}
\label{limit measure}
\lim_{t\to\8} t^\a L_b(t) \E\big[ f(t^{-1}X) \big] = \is f{\L^1}.
\end{equation}
Moreover, the above convergence holds for every bounded function
$f$ such that $0\notin {\rm supp} f$ and $\L^1({\rm Dis}(f))=0$
 (${\rm
Dis}(f)$ is the set of all discontinuities of the function $f$).
In particular
$$
\lim_{t\to\8} t^\a L_b(t) \P\big[|X|>t \big] = \is {{\bf 1}_{\{|\cdot|>1\}}}{\L^1}.
$$
\end{thm}

There are many results describing existence of stationary measures of Markov chains and their tails, especially in the context of general stochastic recursions (see e.g. \cite{DF,G} for one dimensional case and \cite{M1} for multidimensional one). Let us return for a moment to the  example of the autoregressive process \eqref{affine rec}.  It is well-known  that if $\E \log^+\|A_1\|<\8$, then
the Lyapunov exponent
$\lambda=\lim_{n\to\8} \frac 1n \log \|A_1\cdot\ldots\cdot A_n\|$
exists and it is constant a.s. \cite{FK}. Moreover, if $\lambda<0$ and $\E\log^+|B_1|<\8$, then the process $X_n$ converges in distribution to the random vector
\begin{equation}
\label{affine sum}
X = \sum_{n=1}^\8 A_1\cdot\ldots\cdot A_{n-1} B_n,
\end{equation}
 whose law $\nu_1$ is the unique stationary measure of the process $\{X_{1,n}\}$. Properties of the measure $\nu_1$ are well described.
The most significant result is due to Kesten \cite{K}, who proved,
under a number of hypotheses, the main being
   $\lim_{n\to\8}\big(\E\|A_1\cdot\ldots\cdot A_n\|^\a\big)^{\frac 1n}=1$ and $\E |B|^\a <\8$,
 for some $\a>0$, that the  measure $\nu_1$ of $\{X_{1,n}^x\}$  is $\a$-regularly varying at infinity (indeed, Kesten proved weaker convergence, however in this context it turns out to be  equivalent with the definition of $\a$-regularly varying measures, see \cite{Basrak, BL}).
 A short and elegant proof of this result in one dimensional settings was given by Goldie \cite{G}.
   Other multidimensional  results were obtained  in \cite{A,BDGHU,Gui4,KP,LeP}.


However, the theorem above concerns a bit different situation. For
the autoregressive process, Theorem \ref{thm: ogon} deals with the
case when the $B$-part is dominating. If we assume that $B_1$ is
$\a$-regularly varying,
 $\lim_{n\to\8}\big(\E\|A_1\cdot\ldots\cdot A_n\|^\a\big)^{\frac 1n}<1$ (then the Markov chain $X_{1,n}$
is $\a$-geometric) and $\E\|A_1\|^\b<\8$ for some $\b>\a$, then hypotheses of Theorem \ref{thm: ogon} are satisfied and we
conclude that $\nu_1$ is $\a$-regularly varying. In this particular case similar results were
proved in one dimension by Grincevicius \cite{Gri} and Grey \cite{Gr} and in the multivariate
setting in \cite{Roi2} and \cite{RW}. However, \cite{RW} deals with the situation of independent
$A_n$ and $B_n$ and in \cite{Roi2} a particular norm $|\sum _{i=1}^dx_ie_i|=\max_{i=1}^d |x_i|$ is
considered. Theorem \ref{thm: ogon} holds for an arbitrary norm and so it
provides a new result even for the recursion \eqref{affine rec}.

\medskip

Our  approach is more general and it may be applied to a
larger class of Lipschitz recursions. It is valid  for
multidimensional generalizations of the autoregressive process
e.g. for recursions: $X_{2,n} = A_n X_{2,n-1} +B_n +C_n(x)$,
   $X_{3,n} =  \max\{A_n X_{3,n-1}, B_n\}$,
   $X_{4,n} = \max\{A_n X_{4,n-1}, B_n\} +C_n$, where $\max\{x,y\} = (\max\{x_1,y_1\},\ldots,\max\{x_d,y_d\})$, for
   $x,y\in\R^d$. Some of these processes were studied in similar context in one dimension in \cite{G,Gr,M1}. Under
   appropriate assumptions, each of these recursions possesses a unique stationary measure and
   its tail is described by Theorem \ref{thm: ogon}.

\medskip

Let us explain the $\g$-geometricity assumption \eqref{eq: geom1},
which ensures contractivity of the system. The standard approach
to stochastic recursions is to assume that the consecutive random
mappings are contractive in average, i.e. $\E \big[\log {\rm Lip
}(\Phi_n)\big]<0$, where ${\rm Lip }(\Phi_n)$ denote the Lipschitz
coefficient of $\Phi_n$ (see e.g. \cite{DF}). However, in higher
dimensions this approach does not provide sufficiently exact
information. One can easily construct a stochastic recursion where
Lipschitz coefficients of random mappings are larger than one, but
the system still possess some contracitivity properties. For
example, consider on $\R^2$ the autoregressive process, where $A$
is a random diagonal matrix with entries on the diagonal $(2,1/3)$
and $(1/3,2)$ both with probability $1/2$. Then the Lipschitz
coefficient of $A$ is always 2, but since $X_n^x - X_n^y =
A_n\cdot\ldots\cdot A_1(x-y)$, the corresponding Markov chain is $\g$-geometric
for small values of $\g$, thus this is a contracitive system. This
is the reason why to study the autoregressive process in higher
dimensions one has to consider the Lyapunov exponents, not
Lipschitz coefficients. And, this is also the reason, we introduce
in more general settings the concept of $\g$-geometric random
processes.


\medskip

Let $\mu$ be the law  of $A$ and $[\supp\mu]\subseteq {\rm
End}(\R^d)$ be the semigroup generated by the support of $\mu$. It
turns out that in a sense formula \eqref{affine sum} is universal
and, even in the general settings, the tail measures can be
described by  similar expressions. Our next theorem is mainly a
consequence of the previous one, but provides a precise
description of the tail measure $\L^1$. This result is interesting
in its own right, but  will play also a crucial role in the proof
of the limit theorem.

Before stating the theorem let us define a sequence $(\Gamma_n)$ of Radon measures on $\R^d\setminus\{0\}$ as follows. Let $\Gamma_1$ be the tail measure of $\Phi(0, B^1)$ (we will prove in Lemma \ref{app} that $\Phi(0, B^1)$ is $\a$-regularly varying). For $n\ge 2$, we define
$\langle f, \Gamma_n\rangle=\E\big[\langle f\circ A_{2}\circ\ldots\circ A_n, \Gamma_1\rangle\big]$.
\begin{thm}\label{thm: ogon1}
Suppose the assumptions of Theorem \ref{thm: ogon} are satisfied.
%
%
%
If $\Phi(x, 0)=x$ for every $x\in\overline{[\supp\mu]\cdot\Phi[\{0\}\times\supp \Lambda_b]}$, and $\lim_{n\to\8}\left(\E\|A_1\cdot\ldots\cdot A_n\|^{\a}\right)^{\frac 1 n}<1$, then the tail measure $\Lambda^1$ defined in \eqref{limit measure} can be expressed as
 \begin{align}\label{ogon1a}
    \langle f, \L^1\rangle= \sum_{k=1}^\8 \langle f, \Gamma_k\rangle =\langle f, \Gamma_1\rangle+\E\bigg[\sum_{k=2}^{\8} \langle f\circ A_{2}\circ\ldots\circ A_k , \Gamma_1\rangle\bigg].
 \end{align}
Furthermore, the measures $\Gamma_n$ are $\a$-homogeneous and their spherical measures satisfy
\begin{align}\label{ogon1b}
    \E\bigg[\int_{\mathbb{S}^{d-1}}f\left(A\ast\omega\right)
    |A\omega|^{\a}\sigma_{\Gamma_n}(d\omega)\bigg]
    =\int_{\mathbb{S}^{d-1}}f(\omega)\sigma_{\Gamma_{n+1}}(d\omega),
\end{align}
 for every $n\in\N$ and $f\in C(\mathbb{S}^{d-1})$, where $A*\omega=\frac{A\omega}{|A\omega|}$. In particular,  the spherical measure of $\L^1$ is given by
 \begin{align}\label{ogon1c}
    \sigma_{\L^1}(d\omega)=\sum_{n=1}^{\8}\sigma_{\Gamma_n}(d\omega).
\end{align}
\end{thm}

\begin{rem} \label{remt2}
The condition: $\Phi(x, 0)=x$ for every
$x\in\overline{[\supp\mu]\cdot\Phi[\{0\}\times\supp\L_b]}\subseteq\R^d$
is only a technical assumption which can be easily verified in
many cases. Indeed, in the case of the  recursion \eqref{affine
rec}, we know that $\Phi(x, y)=x+y$ and then one has nothing to
check. In the case of the recursion \eqref{max rec}, $\Phi(x,
y)=\max\{x, y\}$ and then $\Phi(x, 0)=x$ holds only for $x\in[0,
\8)$, so we need to know whether
$\overline{[\supp\mu]\cdot\Phi[\{0\}\times\supp\L_b]}\subseteq[0,
\8)$. It is clear that the inclusion depends on the underlying
random variables $A$ and $B^1$, and the sufficient  assumptions
are $\P[A\ge0]=1$ and $\lim_{t\to\8}t^{\a}\P[B^1>t]=c>0$
\end{rem}


   \medskip


In the second part of the paper we study  behavior of the Birkhoff sums
$S_n^x$. We prove that if $\a \in (0,2)$ then there are
constants $d_n, a_n$ such that
$
a_n^{-1}S_n^x-d_n$ converges in law to an $\a $-stable random variable.
In order to state our results
we need some further hypotheses and definitions.

\medskip

The normalization of partial sums will be given by the sequence of numbers $a_n$
defined by the formula
$$a_n=\inf\big\{t>0:
\nu\{x\in\R^d: |x|>t\})\le  1/n\big\},$$ where $\nu$ is the stationary distribution of
$\{X_n^x\}$. One can easily prove that (see Theorem 7.7  in \cite{dur} page 151)
\begin{equation}
\label{an}
    \lim_{n\to\8}n\P(|X|>a_n)=1\ \ \mbox{and}\ \ \ \lim_{n\to\8}\frac{a_n^{\a}L_b(a_n)}{n}= \is {{\bf 1}_{\{|\cdot|>1\}}}{\L^1}=c>0,
\end{equation}
for $\Lambda^1$ being the tail measure of the stationary solution $X$ as in Theorem \ref{thm: ogon}.

\medskip

The characteristic functions of limiting random variables depend on the measure $\L^1$.
However, in their description another Markov chain will play a significant role. Let
$
W_n^x = \ov \Phi_n(W_{n-1}^x)$, where  $W_0^x=x\in \R^d$,
 $\ov\Phi_n (x) = \Phi(A_n x,0)$ and let
$
W(x) = \sum_{k=1}^\8 W_k^x
$.
Then $W_n^x$ is a particular case of recursion \eqref{gen_rec}, with $B_n = 0$. Given $v\in\R^d$
we define $h_v(x) = \E\big[ e^{i\is v{W(x)}} \big]$.

\medskip

Our next result is

\begin{thm}\label{thm: LT}
Suppose that the assumptions of Theorem \ref{thm: ogon} 
are satisfied for some $\a\in(0,2)$. Assume additionally that $\Phi$ is a Lipschitz mapping
and that there is a finite constant $C>0$ such that $|B^2|\le C$ a.e.  Then
the sequence $a_n^{-1}S_n^x-d_n$ converges in law to an $\a$-stable
random variable with the Fourier transform $\ups_\a(tv) = \exp C_\a(tv)$, for
\begin{eqnarray*}
C_\a(tv) &=&
\frac {t^{\a}}{c}
\int_{\R^d}\Big(\big(e^{i\langle v, x\rangle}-1\big)h_v(x)\Big)\L^1(dx), \qquad \mbox { if }\a \in (0, 1);\\
C_1(tv) &=&
\frac {t}{c}
\int_{\R^d}\left(\big(e^{i\langle v, x\rangle}-1\big)h_v(x)-\frac{i\langle v,
x\rangle}{1+|x|^2}\right)\L^1(dx) -\frac{it\log t\langle v, m_{\sigma_{\L^1}}\rangle}{c}, \qquad \mbox { if }\a = 1;\\
C_\a(tv) &=&
\frac {t^{\a}}{c}
\int_{\R^d}\Big(\big(e^{i\langle v, x\rangle}-1\big)h_v(x)-i\langle v,
x\rangle\Big)\L^1(dx), \qquad \mbox { if }\a \in (1,2);
\end{eqnarray*}
where $t>0$, $v\in\mathbb{S}^{d-1}$,  $c$ is the constant defined in \eqref{an} and $m_{\sigma_{\L^1}}=\int_{\mathbb{S}^{d-1}}\o\sigma_{\L^1}(d\o)$ and $\sigma_{\L^1}$ is the spherical measure of the tail measure $\L^1$ defined in Theorem \ref{thm: ogon},
\begin{itemize}
\item if $\a\in (0,1)$, $d_n=0$;
\item if $\a=1$, $d_n = n\xi(a_n^{-1})$,
$\xi(t)=\int_{\R^d}\frac{tx}{1+|tx|^2}\nu(dx)$;
\item if $\a\in(1,2)$, $d_n=a_n^{-1}nm$, for $m=\int_{\R^d}x\nu(dx)$.
\end{itemize}
The functions $C_\a$ satisfy $C_\a(tv) = t^\a C_\a(v)$ for $\a\in(0, 1)\cup(1, 2)$.

Moreover, if $\lim_{n\to\8}\left(\E\|A_1\cdot\ldots\cdot A_n\|^{\a}\right)^{\frac 1 n}<1$, $\Phi(x,0)=x$ for every $x\in\overline{[\supp\mu]\cdot\supp\nu}$, and $\Phi[\{0\}\times\supp\sigma_{\L_b}]$ is not contained in any proper subspace of $\R^d$, then the limit laws
are fully nondegenerate, i.e. $\Re C_\a(tv)<0$ for every $t>0$ and $v\in\mathbb{S}^{d-1}$ and $\a\in(0, 2)$.
\end{thm}
\begin{rem}\label{remt3}
The condition: $\Phi(x,0)=x$ for every
$x\in\overline{[\supp\mu]\cdot\supp\nu}$, requires an explanation
as in Remark \ref{remt2}. It is obvious if $\Phi(x,y)=x+y$. For
instance, if $\Phi(x,y)=\max\{x, y\}$, then $\Phi(x,0)=x$ for $x\in[0,
\8)$, it is
sufficient to assume
$\P[A\ge0]=1$, $\E[A^{\a}]<1$ and $\lim_{t\to\8}t^{\a}\P[B^1>t]=c>0$.
\end{rem}


If $\a>2$ then  $\frac{S_n^x-nm}{\sqrt n}$ converges to a normal law which is a straightforward
application of the martingale method, see \cite{Be, Roi3,WW} and the references given there.
Let us underline that the theorem above concerns dependent  random variables with infinite variance. In the context of stochastic recursions similar problems were studied e.g. in \cite{BJMW, BDG, GL, M1}.
Our proof of  Theorem \ref{thm: LT} is based on the spectral method, introduced by Nagaev in 50's to prove  limit theorems for Markov chains. This method has been strongly developed recently and it has been used in the context of
 limit theorems related to stochastic recursions, see e.g. \cite{BDG, GL, HH, M1}.

\medskip

Throughout the whole paper, unless otherwise stated, we will use the convention that $C>0$ stands for a large positive constant
whose value varies from occurrence to occurrence.

\section*{Acknowledgements}
The authors are grateful to the referees for a very careful reading
of the manuscript and useful remarks that lead to the improvement
of the presentation.

\section{Tails of random recursions}

First we will prove existence and  uniqueness of the stationary measure for the Markov chain
$\{X_n^x\}$ defined in \eqref{gen_rec} as well as some further properties of $\gamma$-geometric Markov
chains that will be used in the sequel.
Following classical ideas, going back to Furstenberg \cite{F} (see also \cite{DF}), we
 consider the backward process $Y_n^x =
\Phi_1\circ\ldots \circ \Phi_n(x)$, which has the same law as $X_n^x$. The process $\{Y_n^x\}$ is not a Markov chain, however
sometimes it is more comfortable to use than $\{X_n^x\}$, e.g. it allows
conveniently to construct the stationary distribution of $\{X_n^x\}$.
Notice that since $X_n^x$ is $\gamma$-geometric, then  $Y_n^x $ is as well, i.e.
\begin{equation}
\label{ylip}
\E\big[  |Y_n^{x} - Y_n^y|^{\gamma} \big] \le C \rho^n|x-y|^{\gamma},\quad x,y\in\R^d,n\in\N,
\end{equation} for $C$ and $\rho$ being as in \eqref{eq: geom1}.

If $x=0$ we  write for simplicity $Y_n$ instead of $Y_n^x$. To emphasize the role of the starting point, which can be sometimes a random variable $X_0$, we  write
$X_n^{X_0}=\Phi_{n}\circ\ldots\circ\Phi_{1}(X_0)$ and $Y_n^{X_0}=\Phi_{1}\circ\ldots\circ\Phi_{n}(X_0)$, where $X_0$ is an arbitrary initial random variable.

\begin{lem}
\label{lem: geom} Let $\{X_n^x\}$ be a Markov chain generated by a
system of random functions, which is $\g$-geometric and satisfies
$\E|X_1|^\d<\8$, for some positive constants $\g,\d>0$. Then there
exists a unique stationary measure $\nu$ of $\{X_n^x\}$  and for any initial random variable $X_0$, the process $\{X_n^{X_0}\}$ converges in distribution to $X$ with the law $\nu$.

Moreover, if additionally $\E|X_0|^\b<\8$ and
$\E\big|X_1^{X_0}\big|^\b<\8$ for some $\b<\g$, then
\begin{equation}
\label{eq: sup}
\sup_{n\in\N} \E|X_n^{X_0}|^\b < \8.
\end{equation}
\end{lem}
\begin{proof}
Take $\eps=\min\{1,\d,\g\}$, then the Markov chain $X_n=X_n^0$ is
$\eps$-geometric. To prove convergence in distribution of $X_n$ it
is sufficient  to show that $Y_n$ converges in $L^{\eps}$.
For this purpose we prove that $\{Y_n\}$ is a Cauchy sequence in
$L^{\eps}$. Fix $n\in\N$, then for any $m>n$ we have
\begin{eqnarray*}
\E\Big[ |Y_m-Y_n|^{\eps} \Big] &\le& \sum_{k=n}^{m-1}\E\Big[ |Y_{k+1}-Y_k|^{\eps} \Big]
= \sum_{k=n}^{m-1}\E\Big[ |Y_{k}^{\Phi_{k+1}(0)}-Y_k|^{\eps} \Big]\\
&\le& C \sum_{k=n}^{m-1} \rho^k \E|\Phi_{k+1}(0)|^{\eps} \le \frac{C\E|X_1|^{\eps} }{1-\rho}\cdot \rho^n.
\end{eqnarray*}
This proves that $Y_n$ converges in $L^{\eps}$, hence also in
distribution, to a random variable $X$. Therefore, $X_n^x$
converges in distribution to the same random variable $X$,
for every $x\in \R^d$.

To prove uniqueness of the stationary measure
assume that there is  another stationary measure $\nu'$. Then, by the
Lebesgue theorem, for every bounded continuous function $f$:
$$
\nu'(f) =  \int_{\R^d}\E \big[ f(X_n^x)   \big]\nu'(dx)
\ _{\overrightarrow{n\to\8}}\ \int_{\R^d} \E\big[f(X)\big]\nu'(dx) = \nu(f),
$$ hence $\nu=\nu'$. The same arguments prove that the sequence $X_n^{Z}$ converges
in distribution to $X$ for any initial random variable $Z$ on $\R^d$.

\medskip

To prove the second part of the lemma, let us consider two cases. 
Assume that $\b<\gamma\le1$, then we write
\begin{eqnarray*}
\E \big|Y_n^{X_0}\big|^\b &\le&  
\sum_{k=0}^{n-1}\E \big|Y_k^{X_0}-Y_{k+1}^{X_0}\big|^\b
+  \E \big|{X_0}\big|^\b\ \le\ \sum_{k=0}^{n-1}  \rho^k   \E \big|X_1^{X_0} - {X_0}\big|^\b
+  \E \big|{X_0}\big|^\b \le C<\8.
\end{eqnarray*}
If $\gamma>1$, it is enough to take $1\le\b<\gamma$ and apply H\"{o}lder inequality, i.e.
\begin{eqnarray*}
\Big( \E \big|Y_n^{X_0}\big|^\b\Big)^{\frac 1\b} &\le& \sum_{k=0}^{n-1} \Big( \E
\big|Y_k^{X_0}-Y_{k+1}^{X_0}\big|^\b\Big)^{\frac 1\b}
+ \Big( \E \big|{X_0}\big|^\b\Big)^{\frac 1\b}\\
&\le& \sum_{k=0}^{n-1}  \rho^k   \Big( \E \big|X_1^{X_0} - {X_0}\big|^\b\Big)^{\frac 1\b}
+ \Big( \E \big|{X_0}\big|^\b\Big)^{\frac 1\b} \le C<\8.
\end{eqnarray*}
\end{proof}

Before we formulate the next lemma, notice that if a random variable $W$ is regularly
varying, then
\begin{equation}
\label{sup}
\sup_{t>0}\Big\{   t^\a L(t) \P\big[|W|>t \big]   \Big\} < \8.
\end{equation}
Moreover, if $L$ is a slowly varying function which is bounded away from zero and infinity on
any compact interval then, by  Potter's Theorem (\cite{BGT}, p. 25), given $\d >0$ there is a finite constant $C>0$
such that
\begin{equation}\label{slow}
\sup _{t>0}\frac{L(t)}{L(\lambda t)}\leq C \max \big\{\l ^{\d }, \l
^{-\d}\big\},
\end{equation}
for every $\l >0$.

\medskip

The following  lemma,  is a multidimensional generalization
of Lemma 2.1 in \cite{DR}.

\begin{lem}\label{app}
Let $Z_1,Z_2\in\R^d$ be $\a$-regularly varying random variables
 with the tail measures $\L_1$, $\L_2$, respectively, (with the same slowly varying function $L_b$ which is bounded away from zero and infinity on
any compact interval), such that
 \begin{equation}\label{niez}
\lim _{t\to \infty} t^\a L_b (t) \P \big[|Z_1| >t, |Z_2| >t\big]=0.
\end{equation}
Then the random variable $(Z_1,Z_2)$ valued in $\R^d\times\R^d$
is regularly varying with index $\a$ and its  tail measure $\L$ is
defined by: $$
   \langle F, \L \rangle=\langle F(\cdot ,0), \L _1\rangle
+\langle F(0, \cdot ), \L _2\rangle, $$
i.e. for every
$F\in C_c((\R ^{d}\times \R^{d})\setminus \{0 \})$:
\begin{equation}
\label{lim}
\lim_{t\to \8} t^\a L_b(t) \E\Big[ F\big(t^{-1}Z_1, t^{-1}Z_2 \big)\Big] = \is F{\L}.
\end{equation}
Moreover, the  formula above is valid for every bounded continuous function $F$
supported outside 0.
\end{lem}
\begin{proof}
 Since every $F\in C_c((\R^{d}\times \R^{d})\setminus \{0 \})$ may be written as a sum of two functions
with supports in $\big(\R^{d}\setminus B_{\eta}(0)\big)\times \R^{d}$
and $\R^{d}\times \big(\R^{d}\setminus B_{\eta}(0)\big)$ respectively, for some  $\eta>0$,
it is enough to consider only one factor of this decomposition.
We assume
that we are in the first case, i.e. supp$F\subseteq \big(\R^{d}\setminus B_{\eta}(0)\big)\times \R^{d}$.
Then to obtain the result for such a function it is enough to justify
that
\begin{equation}
\label{eq:2}
\lim_{t\to\8} t^\a L_b(t)\E\Big[ F\big(t^{-1}Z_1,t^{-1}Z_2\big) -  F\big(t^{-1}Z_1,0\big) \Big] = 0.
\end{equation}
Fix $\eps>0$ and write
\begin{multline*}
t^\a L_b(t)\bigg|\E\Big[ F\big(t^{-1}Z_1,t^{-1}Z_2\big) -  F\big(t^{-1}Z_1,0\big) \Big] \bigg| \\
\le t^\a L_b(t)\E\Big[ \big|F\big(t^{-1}Z_1,t^{-1}Z_2\big) \big| {\bf 1}_{\{|Z_2|>\eps t\}} \Big]
+t^\a L_b(t)\E\Big[ \big|F\big(t^{-1}Z_1,0\big) \big| {\bf 1}_{\{|Z_2|>\eps t\}} \Big]\\
+t^\a L_b(t)\E\Big[ \big|F\big(t^{-1}Z_1,t^{-1}Z_2\big) - F\big(t^{-1}Z_1,0\big) \big| {\bf 1}_{\{|Z_2|\le \eps t\}} \Big]
\end{multline*}
We denote the consecutive expressions in the sum above by $g_1(t),g_2(t),g_3(t)$, respectively.
Taking $\l=\min\{\eta,\eps\}$, by \eqref{slow} and \eqref{niez} we obtain
\begin{multline*}
0\le\lim_{t\to\8} g_1(t) \le \lim_{t\to\8} t^\a L_b(t)\|F\|_\8 \P\big[|Z_1|>\eta t, |Z_2|> \eps t  \big]\\
\le \|F\|_\8\cdot \sup_{t>0} \frac{ L_b(t)}{ L_b(\l t)} \cdot \lim_{t\to\8} \Big(t^\a L_b(\l t)\P\big[|Z_1|>\l t, |Z_2|> \l t  \big]\Big)
 = 0.
\end{multline*}
Arguing in a similar way as above we deduce that $\lim_{t\to\8} g_2(t)=0$.
Finally, to prove that $g_3$ converges to 0, assume first that $F$
is a Lipschitz function with the Lipschitz coefficient ${\rm
Lip}(F)$. Then by \eqref{sup}
\begin{multline*}
g_3(t) \le {\rm Lip}(F) t^\a L_b(t) \E\Big[ |t^{-1}Z_2| {\bf 1}_{\{ |t^{-1}Z_1|>\eta\}}{\bf 1}_{\{ |t^{-1}Z_2|\le \eps\}} \Big]\\
\le \eps\cdot {\rm Lip}(F) \sup_{t>0}\Big\{  t^\a L_b(t) \P[|t^{-1}Z_1|>\eta] \Big\}
\le C\eps.
\end{multline*}
Passing with $\eps$ to 0, we obtain \eqref{eq:2} for Lipschitz functions.

\medskip

To prove the result for arbitrary functions,
notice first that
\eqref{sup} implies
$$
\sup_{t>0}\Big\{ t^\a L_b (t)\P\big[\eta t < |Z_1| + | Z_2| <Mt\big]\Big\}<\infty.$$
 Now we approximate $F\in C_c\big(\big(\R^{d}\setminus B_{\eta}(0)\big)\times \R^{d}\big)$ by a
 Lipschitz function $G \in C_c\big(\big(\R^{d}\setminus B_{\eta}(0)\big)\times \R^{d}\big)$
such that $\| F-G\|_\8<\eps $. Then
\begin{multline*}
 t^\a L_b (t) \Big | \E \big[ F(t^{-1}Z_1,t^{-1}Z_2 )- F(t^{-1}Z_1,0)\big ]\Big |
 \le
t^\a L_b (t)  \E \big[\big | F(t^{-1}Z_1,t^{-1}Z_2 )- G(t^{-1}Z_1,t^{-1}Z_2)\big|\big ]\\
+t^\a L_b (t) \Big | \E \big[ G(t^{-1}Z_1,t^{-1}Z_2 )- G(t^{-1}Z_1,0)\big ]\Big |
+t^\a L_b (t) \E \big[\big |  F(t^{-1}Z_1,0 )- G(t^{-1}Z_1,0)\big ]\Big |\\
 \le
\eps t^\a L_b (t) \P\big[\eta t < |Z_1| + | Z_2| <Mt\big]
+t^\a L_b (t) \Big | \E \big[ G(t^{-1}Z_1,t^{-1}Z_2 )- G(t^{-1}Z_1,0)\big ]\Big |\\
+\eps t^\a L_b (t) \P\big[\eta t < |Z_1|<Mt\big],
\end{multline*}
hence passing with $t$ to infinity and then with $\eps$ to zero we obtain \eqref{eq:2} and so also \eqref{lim}.

\medskip

To prove the second part of the lemma,
let $F$ be an arbitrary bounded continuous function on $\R^d\times\R^d$ supported outside 0. Assume $\|F\|_\8=1$.
Take $r>0$ and let $\phi_1,\phi_2$ be nonzero functions on $\R^d\times\R^d$ such that $\phi_1+\phi_2=1$,
${\rm supp}\phi_1 \subseteq B_{2r}(0)$ and
${\rm supp}\phi_2 \subseteq B_{r}(0)^c$. Then by \eqref{sup} and \eqref{slow}
\begin{multline*}
\lim_{r\to\8}\sup_{t>0}t^\a L_b(t) \E\big[ (\phi_2 F)(t^{-1}Z_1,t^{-1}Z_2) \big] \le
\lim_{r\to\8}\sup_{t>0} t^\a L_b(t) \Big(\P\big[ |Z_1|>rt \big] + \P\big[ |Z_2|>rt \big]\Big) \\
\le \lim_{r\to\8}\sup_{t>0} r^{-\a}\frac{L_b(t)}{ L_b(rt)} (rt)^\a L_b(rt) \Big(\P\big[ |Z_1|>rt \big] + \P\big[ |Z_2|>rt \big]\Big)
= 0.
\end{multline*}
By \eqref{lim}
$$
\lim_{t\to\8}  t^\a L_b(t) \E\big[ (\phi_1 F)(t^{-1}Z_1,t^{-1}Z_2) \big] = \is {\phi_1 F}{\L}.
$$
Therefore, passing with $r$ to infinity, we obtain \eqref{lim}
for non-compactly supported functions $F$.
\end{proof}

The next lemma when considered for the one dimensional recursion
\eqref{affine rec} is known as Breiman's lemma \cite{B}. In the
multidimensional affine settings the lemma was proved in
\cite{Roi2} (Lemma 2.1). Here we write it in the generality
corresponding to our framework and, at the same time, we present a
simpler proof than in \cite{Roi2}.

\begin{lem}\label{spectral}
Assume that
\begin{itemize}
\item random variables $(A,B)$ and $X\in \R^d$ are independent;
 \item $X$ and $B^1$ are $\a $-regularly varying
with the tail measures  $\L$, $\L_b$, respectively, (with the same slowly varying function $L_b$ which is bounded away from zero and infinity on
any compact interval);
\item $\E \| A\| ^{\b}<\8$ for some $\b
>\a$;
\item there is $\eps_0>0$ such that $\E\big[ (B^3)^{\frac{\a}{\d_0}+\eps_0} \big]<\8$, if $0<\d_0<1$ and  $\E\big[ (B^3)^{\a+\eps_0} \big]<\8$, if $\d_0=0$.
\end{itemize}
Then both  $AX$ and $\Phi(AX,B(X))$ are
 $\a$-regularly varying with the tail measures $\wt \L$ and $\L_1$ respectively, where  $\is f{\wt \L} =  \E\big[ \langle f\circ A,
\Lambda \rangle\big]$ and
\begin{equation}
\label{eq:lambda1}
\is f{\L_1} = \is{f\circ\Phi(\cdot,0)}{\wt\L} + \is{f\circ\Phi(0,\cdot)}{\L_b}.
\end{equation}
\end{lem}
\begin{proof}
First, conditioning on $A$, we will prove that for any bounded function $f$ supported in $\R^d\setminus B_{\eta}(0)$ for some $\eta>0$,
there exists a function $g$ such that
\begin{equation}
\label{eq:1} \sup_{t>0}\Big\{ t^{\a }L_b(t)  \E\big[ f\big(t^{-1}AX\big)|A \big]\Big\} \le
g(A),\quad\mbox{ and }\quad \E[g(A)] <\8.
\end{equation}

Observe that $\sup _{t>0}t^{\a }L_b(t)\P \big[ | X| >t\big]=C<\8$ and assume that
$\mbox{supp}f\subseteq \R ^d\setminus B_{\eta}(0)$, $\eta<1$, and fix $\d<\b-\a$.
If $\| A\| \leq 1$  then, by \eqref{slow}, for every $t>0$
$$
 t^{\a}L_b(t)\E\big[ f\big(t^{-1}AX\big)|A \big]
 \leq \| f \|_\8
 t^{\a }L_b(t)\P \big[ | X| >t\eta \big]\leq C\eta ^{-\a -\d }\|
f\|_\8 =C_1<\8.$$
If $2^n\leq \| A\| \leq 2^{n+1}$ for $n\in\N$ then, again by \eqref{slow}, for every $t>0$
\begin{eqnarray*}
 t^{\a}L_b(t) \E\big[ f\big(t^{-1}AX\big)|A \big]&\leq&
 \| f \|_\8 t^{\a }L_b(t)\P \big[ 2^{n+1}| X|> t\eta \big]\\ &\leq& C2^{(n+1)(\a
+\d )}\eta ^{-\a -\d }\| f\|_\8  = C_2 2^{n(\a+\d )}.
\end{eqnarray*}
Finally, notice that
\begin{eqnarray*}
\E [ g(A)] &\le& C_1  \P[\|A\|\le 1] + C_2\sum_{n=1}^\8 2^{n(\a+\d )} \P \big[ \| A\| \geq 2^{n}\big] \\
&\leq&
C_1  + C_2  \E \| A\| ^{\b}\cdot\sum_{n=1}^\8 2^{n(\a+\d-\b )} <\8,
\end{eqnarray*}
 and the
proof of \eqref{eq:1} is completed. Now in view of \eqref{eq:1} we can easily  prove that $AX$ is regularly varying with index $\a$. Indeed, taking $f\in C_c(\R^d\setminus B_{\eta}(0))$,  conditioning on $A$, and
using dominated convergence theorem we have
$$\lim_{t\to\8} t^\a L_b(t) \E \big[ f(t^{-1} AX) \big] =
\E\bigg[ \lim_{t\to\8} t^\a L_b(t) \E\Big[ (f\circ A)(t^{-1}X)\big| A\Big]\bigg]
= \E\big[\is{f\circ A}{\L}\big] =\is f{\wt\L},
$$
hence $AX$ is $\a$-regularly varying as desired.

\medskip

For the second part of the lemma, we are going to apply Lemma \ref{app},
with $Z_1 = AX$, $Z_2=B(X)$ and the function $f\circ\Phi$. Notice, that since $\Phi(0,0)=0$
the function $f\circ\Phi$ is supported outside 0. It may happen (e.g. when $\Phi(x,y)=x+y$)
that $f\circ\Phi$ is not compactly supported, however it is still a bounded function.
Therefore, we have to prove that $B(X)$ is $\a$-regularly varying with
the tail measure $\L_b$ and \eqref{niez} is satisfied, i.e.
\begin{equation}
\label{eq:11}
\lim_{t\to\8} t^\a L_b(t) \P\big[ |AX|> t, |B(X)|>t
\big] =0.
\end{equation}
To prove that $B(X)$ is $\a$-regularly varying notice that from the first
part of the lemma with $B^3$ instead of $A$ we know that
if $\d_0>0$, then $(B^3)^{\frac1 {\d_0}}X$ is $\a$-regular. Therefore,
$$
\lim_{t\to\8} t^\a L_b(t)
 \P\big[ B^2(X)>t \big]
\le \lim_{t\to\8} t^\a L_b(t)
 \P\big[ (B^3)^{\frac1 {\d_0}}|X|>t ^{\frac 1{\d_0}}\big] = 0,
$$
so $B^2(X)$ is $\a$-regularly varying with the tail measure $0$. If $\d_0=0$, then  $\lim_{t\to\8} t^\a L_b(t)
 \P\big[ B^2(X)>t \big]=0$ can be easily established. Hence applying Lemma \ref{app} for $Z_1=B_1$, $Z_2 = B^2(X)$ and
$f\circ\widetilde{\Phi}$, where $\widetilde{\Phi}(x,y)=x+y$ we deduce
\begin{multline*}
\lim_{t\to\8} t^\a L_b(t) \E\big[ f(t^{-1}B(X))\big]
=  \lim_{t\to\8} t^\a L_b(t) \E\Big[ (f\circ\widetilde{\Phi})\big(t^{-1}B^1,t^{-1}B^2(X)\big)\Big]\\
= \is{(f\circ\widetilde{\Phi})(\cdot,0)}{\L_b}  + \is{(f\circ\widetilde{\Phi})(0,\cdot)}{0} = \is{f}{\L_b}.
    \end{multline*}
In order to prove \eqref{eq:11} take $f(x)=\1_{\{| \cdot |>1\}}(x)$, then applying \eqref{eq:1} and
conditioning on $(A,B^1)$ we obtain
\begin{multline*}
t^\a L_b(t)\P\big[|AX|>t,|B(X)|>t\big] \le
 t^\a L_b(t) \E\big[ f(t^{-1}AX)\1_{\{|B^1|>t/2\}} \big] + t^\a L_b(t)\P\big[|B^2(X)|>t/2\big]\\
\le  \E\Big[ \1_{\{|B^1|>t/2\}} \cdot  \sup_{t>0} t^\a L_b(t) \E\big[  f(t^{-1}AX) |(A,B^1)\big]\Big]
+ t^\a L_b(t)\P\big[|B^2(X)|>t/2\big]\\
\le \E\big[ \1_{\{ |B^1|>t/2 \}} g(A) \big]+ t^\a L_b(t)\P\big[|B^2(X)|>t/2\big].
\end{multline*}
The last expression converges to 0 as $t$ goes to infinity.
Finally, from Lemma \ref{app} we obtain that $\Phi(A,B)(X)$ is
$\a$-regular:
$$\lim_{t\to\8} t^\a L_b(t) \E\Big[ f\big(t^{-1} \Phi(A,B)(X)\big)  \Big]
= \lim_{t\to\8} t^\a L_b(t) \E\Big[ \big(f\circ\Phi\big)\big(t^{-1} AX, t^{-1}B(X)\big)  \Big]
= \is f{\L_1}.
$$
This proves \eqref{eq:lambda1} and completes the proof of the lemma.
\end{proof}

\begin{proof}[Proof of Theorem \ref{thm: ogon}]
Since the stationary solution $X$ does not depend on the choice of the initial random variable $X_0$, without any loss
of generality, we may assume that $X_0$ is $\a$-regularly varying with some
nonzero tail
measure $\L_0$. Then by Lemma \ref{spectral}, for every $n\in\N$, $X_n^{X_0}$ is $\a$-regularly
varying  with the tail measure $\L_n$ satisfying \eqref{eq:lambda1} with $\wt \L_{n-1}$, being the tail measure of $A_n X_{n-1}^{X_0}$,
instead of $\wt\L$.  So, we have to prove that $\L_n$ converges weakly to some measure $\L^1$, which
we can identify as the tail measure of $X$. This measure will be nonzero, since for every $n\in\N$ and positive $f$:
$\is f{\L_n}\ge \is{f\circ\Phi(0,\cdot)}{\L_b}$.
 From now
we will consider the backward process $\{Y_n^x\}$. We may assume that $\d>0$ in \eqref{slow}
is sufficiently small, i.e. $\d<\min\{\a,\g-\a\}$.
Suppose first that $f$ is an $\eps$-H\"older function for $0<\eps<\d$
 and $\mbox{supp} f \subseteq \R^d\setminus B_{\eta}(0)$. By \eqref{ylip} there exist
constants $0<C_0<\8$ and $0<\rho_0<1$ such that
\begin{equation}
\label{eq:unif}
 \E\big[\big| Y_n^{x} - Y_n^{y}  \big|^s\big] \le C_0\rho_0^n |x-y|^s   \quad\mbox{for } s\in\{\g,\a-\d,\a+\d\},
  n\in\N, \mbox{ and } x,y\in \R^d.
 \end{equation}
We will  prove that there are constants $0<C<\8$ and $0<\rho<1$
such that for every $m>n$
\begin{equation}\label{cauchy}
\sup_{t>0}\Big\{ t^{\a}L_b(t)\E \big|f(t^{-1}Y_m^{X_0})-f(t^{-1}Y_n^{X_0})\big|\Big\} \leq C \rho^n.
\end{equation}
We begin by showing that
\begin{equation}
\label{eq:yn}
\sup_{t>0} \Big\{ t^\a L_b(t) \E \big|f\big( t^{-1}Y_k^{X_0} \big) - f\big( t^{-1}Y_k \big)\big|  \Big\} \le C \rho^k,
\end{equation}
for $k\in\N$. We have
\begin{multline*}
\E \big[f(t^{-1}Y_{k}^{X_0})- f(t^{-1}Y_k)\big]
= \E \Big[ \big(f(t^{-1}Y_{k}^{X_0})-f(t^{-1}Y_k)\big) \1_{\{|t^{-1}Y_k| >\frac{\eta}{2}\}}\Big] \\
+ \E \Big[\big (f(t^{-1}Y_{k}^{X_0})- f(t^{-1}Y_k)\big )  \1_{\{|t^{-1}Y_k^{X_0}|>\eta\}} \1
_{\{|t^{-1}Y_k|<\frac{\eta}{2}\}}\Big]=I_1+I_2.
\end{multline*}
Notice that $\E|\Phi_1(0)|^\b<\8$ for every $\b<\a$, hence by \eqref{eq: sup}: $\sup_{k\in\N}\E|Y_k|^\b\le C<\8$.
 Therefore, on the one hand, we have an estimate for small $t>0$
 \begin{eqnarray*}
t^\a L_b(t)|I_1|
&\le& C t^{\a-\eps} L_b(t) \E\bigg[\E\Big[ \big| Y_k^{X_0} - Y_k \big|^{\eps} {\bf 1}_{\{|Y_k|>t\eta/2\}} \Big| X_0\Big]\bigg]\\
&\le& C t^{\a-\eps} L_b(t) \E\big[|X_0|^{\eps}\big]\rho_0^k.
\end{eqnarray*}
 On the other hand, by the H\"older inequality with $p=\frac{\g}{\eps}$, $q=\frac{\g}{\g-\eps}$, conditioning on
 $X_0$ we have an estimate for  sufficiently large $t>0$
\begin{eqnarray*}
t^\a L_b(t)|I_1|
&\le& C t^{\a-\eps} L_b(t) \E\bigg[\E\Big[ \big| Y_k^{X_0} - Y_k \big|^{\eps} {\bf 1}_{\{|Y_k|>t\eta/2\}} \Big| X_0\Big]\bigg]\\
&\le& C t^{\a-\eps} L_b(t) \E\bigg[\E\Big[ \big| Y_k^{X_0} - Y_k \big|^{p\eps}  \Big| X_0\Big] ^{\frac 1p}
\E\Big[{\bf 1}_{\{|Y_k|>t\eta/2\}} \Big| X_0\Big]^{\frac 1q}
 \bigg]\\
&\le& C t^{\a-\eps} L_b(t) \E\bigg[\E\Big[ \big| Y_k^{X_0} - Y_k \big|^{\gamma}  \big| X_0\Big] ^{\frac 1p} \bigg]
\P\big[ |Y_k|>t\eta/2 \big]^{\frac 1q}\\
&\le&  C t^{\a-\eps} L_b(t) \rho_0^{\frac kp} \E|X_0|^{\eps} \cdot t^{-\big( \a - \frac{\eps \d}{\g-\eps} \big)\frac 1q} \E\Big[|Y_k|^{ \a - \frac{\eps \d}{\g-\eps}}\Big]^{\frac 1 q}\\
&\le&  C L_b(t) t^{\frac {1}p(\a+\d-\g)}  \rho_0^{\frac kp}.
\end{eqnarray*}
Finally, we have obtained
\begin{eqnarray*}
t^\a L_b(t)|I_1|
&\le& C L_b(t)\min\Big\{  t^{\a-\eps},
 t^{\frac {1}p(\a+\d-\g)}  \Big\}\rho_0^{\frac kp}.
\end{eqnarray*}
Denote by $\wt L_n$ the Lipschitz coefficient of $\Phi_1\circ\cdots\circ \Phi_n$.
Since $X_0$ is $\a$-regularly varying, by \eqref{sup} and \eqref{slow} we obtain
\begin{eqnarray*}
t^\a L_b(t) |I_2| &\le& 2\|f\|_\8 t^\a L_b(t) \P\Big[ \big| Y_k^{X_0} - Y_k \big| > t\eta/2 \Big]\\
&\le& 2\|f\|_\8 t^\a L_b(t) \P\Big[ \wt L_k |X_0| > t\eta/2 \Big]\\
&\le& C\|f\|_\8 \E\Bigg[  \wt L_k^\a \frac{L_b(t)}{L_b\big(\frac{t\eta}{2\wt L_k}\big)} \E\bigg[
\bigg(\frac {t\eta}{2\wt L_k}\bigg)^\a L_b\bigg(\frac {t\eta}{2\wt L_k}\bigg)
{\bf 1}_{\big\{|X_0|> \frac{t\eta}{2\wt L_k}\big\}}\Big| \wt L_k\bigg] \Bigg]\\
&\le& C\|f\|_\8 \E\Big[ \wt L_k^{\a+\d} + \wt L_k^{\a-\d}\Big] \le C\|f\|_\8 \rho_0^k.
\end{eqnarray*}
Hence, we deduce \eqref{eq:yn} and in order to prove
\eqref{cauchy} it is enough to justify
\begin{equation}\label{cauchyy}
\sup_{t>0} \Big\{ t^{\a}L_b(t)\E\Big[ \big|f(t^{-1}Y_m)-f(t^{-1}Y_n)\big|\Big]\Big\} \leq C \rho^n, \qquad m>n.
\end{equation}
For this purpose we decompose
$$
f(t^{-1}Y_m)-f(t^{-1}Y_n)=\sum_{k=n}^{m-1}\big(f(t^{-1}Y_{k+1})-f(t^{-1}Y_k)\big),$$
and next we estimate $\E[f(t^{-1}Y_{k+1})-f(t^{-1}Y_k)]$
using exactly the same arguments as in \eqref{eq:yn}, with $Y_{k+1}=Y_k\circ\Phi_{k+1}$ instead of $Y_k^{X_0}$
and $\Phi_{k+1}(0)$ instead of $X_0$. Thus we obtain that
\begin{equation}
\label{eq:yk}
\sup_{t>0} \Big\{ t^\a L_b(t) \E \big|f\big( t^{-1}Y_{k+1} \big) - f\big( t^{-1}Y_k \big)\big|  \Big\} \le C \rho^k,
\end{equation}
which in turn implies \eqref{cauchyy} and hence \eqref{cauchy}.
Now letting $m\to \8$ we have
\begin{equation}
\label{eq:yx}
\sup_{t>0}\Big\{
t^{\a}L_b(t)\E\big[ |f(t^{-1}X)-f(t^{-1}Y_n^{X_0})| \big] \Big\}\leq
C \rho^n.
\end{equation}
We know that, for every $n\in\N$, $Y_n^{X_0}$ is
$\a$-regularly varying with the tail measure $\L_n$. Moreover,
in view of \eqref{cauchy}, the sequence $\L_n(f)$ is a Cauchy
sequence, hence it converges. Let $\L^1(f)$ denotes the limit of $\L_n(f)$. In view of \eqref{eq:yx}, for every $n\in\N$,
we have
\begin{multline*}
\limsup_{t\to\8} \Big| t^\a L_b(f) \E\big[f(t^{-1} X) \big]  - \L^1(f) \Big|
\le \limsup_{t\to\8}  t^\a L_b(f) \E\big[|f(t^{-1} X) - f(t^{-1}Y_n^{X_0})|  \big] \\
+ \lim_{t\to\8} \Big| t^\a L_b(f) \E\big[f(t^{-1} Y_n^{X_0}) \big] - \L_n(f) \Big|
+ \big| \L_n(f) -\L^1(f) \big| \le C\rho^n+\big| \L_n(f) -\L^1(f) \big|,
\end{multline*}
and so letting $n\to \8$
\begin{equation} \label{limit}
  \lim_{t\to\8}  t^\a L_b(f) \E\big[f(t^{-1} X)  \big] = \L^1(f),
  \end{equation}
for any $\eps$-H\"older function.

Finally, take a continuous function $f$ compactly supported in $\R^d\setminus
B_{\eta}(0)$ for some $\eta>0$, and fix $\d>0$. Then there exists an $\eps$-H\"older function $g$ supported in
$\R^d\setminus B_{\eta}(0)$ such that $\|f-g\|_\8\le \d$. Moreover, let $h$ be an $\eps$-H\"older function, supported in $\R^d\setminus B_{\eta/2}(0)$,
 such that $\d h\ge  |f-g|$. To define $\L^1(f)$ we
will first prove an inequality similar to \eqref{cauchy}. Notice
that
\begin{multline*}
\sup_{t>0} \Big\{ t^{\a}L_b(t)\E \big|f(t^{-1}Y_m)-f(t^{-1}Y_n)\big|\Big\}
\le
\sup_{t>0} \Big\{ t^{\a}L_b(t)\E \big|f(t^{-1}Y_m)-g(t^{-1}Y_m)\big|\Big\} \\
+\sup_{t>0} \Big\{ t^{\a}L_b(t)\E \big|g(t^{-1}Y_m)-g(t^{-1}Y_n)\big|\Big\}
+\sup_{t>0} \Big\{ t^{\a}L_b(t)\E \big|g(t^{-1}Y_n)-f(t^{-1}Y_n)\big|\Big\} \\
\le \d \L_m(h)  + C\rho^n  +  \d \L_n(h),
\end{multline*}
hence $\L_n(f)$ is a Cauchy sequence, since $\d>0$ is arbitrary. Denote its limit by
$\L^1(f)$. Then $\L^1$ is a well defined Radon measure on
$\R^d\setminus\{0\}$.

\medskip

To prove the second part of the theorem we proceed as at the
end of the proof of Lemma \ref{app}, obtaining \eqref{limit} for
bounded continuous functions supported outside $0$.  By the
Portmanteau theorem we have also  \eqref{limit} for every bounded
function $f$ supported outside 0 and such that $\L^1({\rm
Dis}(f))=0$. Finally, since $\L^1$ is $\a$-homogeneous, it can be
written in the form \eqref{radial}, hence we have
 $\L^1\big({\rm Dis}({\bf 1}_{\{|\cdot|>1\}})\big)=0$, and the proof of Theorem \ref{thm: ogon} is completed.
\end{proof}

\begin{proof}[Proof of Theorem \ref{thm: ogon1}] Since the stationary solution $X$ does not depend on the choice of the starting point we may assume, without any loss of generality, that $X_0=0$ a.s., then in view of Lemma \ref{spectral} we know that $X_1=\Phi(A_1X_0, B_1(X_0))=\Phi(0, B_1^1)$ is $\alpha$-regularly varying with the tail measure $\Lambda_1$ (notice $\Lambda_1 = \Gamma_1$). Applying Lemma \ref{spectral}, to the random variable $X_2=\Phi(A_2X_1, B_2(X_1))$, we can express its tail measure $\Lambda_2$ in the terms of $\Lambda_1$. Indeed,
\begin{align*}
    \langle f, \L_2\rangle&=\langle f\circ\Phi(\cdot, 0), \widetilde{\L}_1\rangle
    +\langle f\circ\Phi(0, \cdot), \L_b\rangle\\
    &=\E\big[ \langle f\circ\Phi(A_2( \cdot ), 0) , \L_1\rangle\big]+\langle f, \L_1\rangle
   = \E\big[ \langle f\circ A_2 , \L_1\rangle\big]+\langle f, \L_1\rangle,
\end{align*}
since $\Phi(x, 0)=x$ for every $x\in\overline{[\supp\mu]\cdot\Phi[\{0\}\times\supp\L_b]}\subseteq\R^d$ and by the definition $\langle f\circ\Phi(0, \cdot), \L_b\rangle=\langle f, \L_1\rangle$.
 If $\L_n$ denotes the tail measure of $X_n$, then an easy induction argument proves  $$\langle f, \L_n\rangle=\E\bigg[\sum_{k=2}^n\langle f\circ A_n\circ\ldots\circ A_k, \L_1 \rangle\bigg]+\langle f, \L_1\rangle, \quad n\in \N.$$
To prove \eqref{ogon1a}, notice that since
$X_n$ has the same law as $Y_n$ and hence
$$\E\bigg[\sum_{k=2}^n\langle f\circ A_n\circ\ldots\circ A_k, \L_1 \rangle\bigg]=\E\bigg[\sum_{k=2}^n\langle f\circ A_2\circ\ldots\circ A_k, \L_1 \rangle\bigg] = \E\bigg[\sum_{k=2}^n\langle f, \Gamma_k \rangle\bigg], $$ for every $n\in\N$. Therefore, we have
 \begin{align}
  \nonumber t^{\a}L_b(t)\E\big[f(t^{-1}X)\big]  -\left(\langle f, \Gamma_1\rangle+\E\bigg[\sum_{k=2}^{\8} \langle f , \Gamma_k\rangle\bigg]\right)
    =t^{\a}L_b(t)\E\big[f(t^{-1}X)\big]-t^{\a}L_b(t)\E\big[f(t^{-1}Y_n)\big]\\
    \label{1star} +t^{\a}L_b(t)\E\big[f(t^{-1}X_n)\big]-\left(\langle f, \Gamma_1\rangle+\E\bigg[\sum_{k=2}^n\langle f, \Gamma_k \rangle\bigg]\right)\\
   \nonumber+\E\bigg[\sum_{k=n+1}^{\8} \langle f, \Gamma_k\rangle\bigg].
 \end{align}
By  \eqref{eq:yx} there exist constants $0<C<\8$ and $0<\rho<1$ such that for every $n\in\N$
 \begin{equation}
 \label{2star}
    \sup_{t>0}\Big|t^{\a}L_b(t)\E\big[f(t^{-1}X)\big]-t^{\a}L_b(t)\E\big[f(t^{-1}Y_n)\big]\Big|\le C\rho^n.
 \end{equation}
 Reasoning as in the first part of the proof of Theorem \ref{thm: ogon} one can prove that for every  $\varepsilon>0$
 there is $t_{\varepsilon}>0$ such that for every $t\ge t_{\varepsilon}$
 \begin{equation}
 \label{3star}
    \left|t^{\a}L_b(t)\E\big[f(t^{-1}X_n)\big]-\left(\langle f, \Gamma_1\rangle+\E\bigg[\sum_{k=2}^n\langle f, \Gamma_k \rangle\bigg]\right)\right|< \varepsilon.
 \end{equation}
Finally assume that $\supp f\subseteq \R^d\setminus B_{\eta}(0)$ for some $\eta>0$, then
\begin{equation}
\label{4star}
\begin{split}
    \left|\E\bigg[\sum_{k=n+1}^{\8} \langle f, \Gamma_k\rangle\bigg]\right|&\le\|f\|_{\8}
    \E\bigg[\sum_{k=n+1}^{\8}\int_{\R^d\setminus\{0\}}\mathbf{1}_{\{y\in\R^d: |y|>\eta\|A_{2}\circ\ldots\circ A_k\|^{-1}\}}(x)\Gamma_1(dx)\bigg]\\
    &\le\eta^{-\a}\|f\|_{\8}\E\bigg[\sum_{k=n+1}^{\8}\|A_{2}\circ\ldots\circ A_k\|^{\a}\big]\bigg]\ _{\overrightarrow{n\to\8}}\ 0,
    \end{split}
\end{equation}
since $\lim_{n\to\8}\left(\E\|A_1\circ\ldots\circ A_n\|^{\a}\right)^{\frac 1 n}<1$.
Combining \eqref{1star} with  \eqref{2star}, \eqref{3star} and \eqref{4star} we obtain \eqref{ogon1a}.

Now take $f\in C_c(\R^d\setminus\{0\})$  of the form $f(r\omega)=f_1(r)f_2(\omega)$, where $r>0$, $\omega\in\mathbb{S}^{d-1}$, $f_1\in C_c((0, \8))$ and $f_2\in C(\mathbb{S}^{d-1})$. In view of Lemma \ref{spectral} we obtain
\begin{eqnarray*}
    \left\langle f_1, \frac{dr}{r^{\a+1}}\right\rangle\langle f_2, \sigma_{\Gamma_n}\rangle&=&\langle f, \Gamma_n\rangle =
    \E\bigg[\int_{\R^d\setminus\{0\}}f(A_2\circ\ldots\circ A_nx)\Gamma_1(dx)\bigg]\\
    &=&\E\bigg[\int_{0}^{\8}\int_{\mathbb{S}^{d-1}}
    f_1(|A_2\circ\ldots\circ A_n\omega|r)f_2((A_2\circ\ldots\circ A_n)*\omega)\sigma_{\Gamma_1}(d\omega)\frac{dr}{r^{\a+1}}\bigg]\\
    &=&\left\langle f_1, \frac{dr}{r^{\a+1}}\right\rangle
    \E\bigg[\int_{\mathbb{S}^{d-1}}
    |A_2\circ\ldots\circ A_n\omega|^{\a}f_2((A_2\circ\ldots\circ A_n)*\omega)\sigma_{\Gamma_1}(d\omega)\bigg],
\end{eqnarray*}
where $A*\omega=\frac{A\omega}{|A\omega|}$ hence we have proved
\begin{align*}
   \langle f_2, \sigma_{\Gamma_n}\rangle=
    \E\bigg[\int_{\mathbb{S}^{d-1}}
    |A_2\circ\ldots\circ A_n\omega|^{\a}f_2\left((A_2\circ\ldots\circ A_n)*\omega\right)\sigma_{\Gamma_1}(d\omega)\bigg],
\end{align*}
Finally to prove  \eqref{ogon1b} we write
\begin{multline*}
    \E\bigg[\int_{\mathbb{S}^{d-1}}f\left(A*\omega\right)
    |A\omega|^{\a}\sigma_{\Gamma_n}(d\omega)\bigg]\\
    =\E\bigg[\int_{\mathbb{S}^{d-1}}f\left(A*((A_2\circ\ldots\circ A_n)*\omega)\right)
    |A((A_2\circ\ldots\circ A_n)*\omega)|^{\a}|A_2\circ\ldots\circ A_n\omega|^{\a}\sigma_{\Gamma_1}(d\omega)\bigg]\\
    =\E\bigg[\int_{\mathbb{S}^{d-1}}f\left((A_2\circ\ldots\circ A_{n+1})*\omega)\right)
    |A_2\circ\ldots\circ A_{n+1}\omega|^{\a}\sigma_{\Gamma_1}(d\omega)\bigg]
    =\int_{\mathbb{S}^{d-1}}f(\omega)\sigma_{\Gamma_{n+1}}(d\omega).
\end{multline*}
Formula \eqref{ogon1c} is a simple consequence of \eqref{ogon1a} and the calculations stated above. This completes the proof of Theorem \ref{thm: ogon1}.

\end{proof}

\section{The Limit Theorem}

Let $\mathcal{C}(\R^d)$ be the space of continuous functions on $\R^d$.  Given positive parameters $\rel$ we introduce two
Banach spaces $\ct(\R^d)$ and $\bl1(\R^d)$ defined as follows

\begin{align*}
\ct=\ct(\R^d)&=\left\{f\in\mathcal{C}(\R^d): |f|_{\rho}=\sup_{x\in\R^d}\frac{|f(x)|}{(1+|x|)^{\rho}}<\8\right\},\\
\bl1=\bl1(\R^d)&=\{f\in\mathcal{C}(\R^d): \|f\|_{\rel}=|f|_{\rho}+[f]_{\el}<\8\},
\end{align*}
where
\begin{align*}
[f]_{\el}=\sup_{x\not= y}\frac{|f(x)-f(y)|}{|x-y|^{\epsilon}(1+|x|)^{\l}(1+|y|)^{\l}}.
\end{align*}

\medskip

On $\ct$ and $\bl1$ we consider the Markov operator
$Pf(x)=\E\big[f(X_1^x)\big]$ and its Fourier
perturbations
\begin{align*}
    P_{t, v}f(x)=\E\left[e^{it\is v{X_1^x}}f(X_1^x)\right],
\end{align*}
where $x\in\R^d$, $v\in\Sd$ and $t>0$. Notice that $P_{0, v}=P$.
The operators will play a crucial role in the proof, since one can
prove by induction that
$$
P_{t,v}^n f(x) = \E\Big[ e^{it\is v{S_n^x}} f(X_n^x)  \Big].
$$
So, the characteristic function of $a_n^{-1}S_n - d_n$ is just
$$
\E\Big[ e^{it\is v{a_n^{-1}S_n - d_n}}  \Big] = P^n_{t a_n^{-1},v}{\bf 1}(x) e^{-it\is v{d_n}}.
$$
Therefore, to prove the theorem one has to consider $P_{t,v}^n$ for large $n$ and small $t$,
what reduces the problem to  describing   spectral properties of the operators $P_{t,v}$
on the Banach space $\bl1$.

\medskip

Next we define another family of Fourier operators
$$
T_{t,v} f(x)=\Dt^{-1}P_{t, v}\Dt f(x), \qquad t>0,
$$ where
 $\Dt$ is  the dilatation
operator defined by   $\Dt f(x)=f(tx)$.
This family is related to the dilated Markov chain $\{X_{n,t}^x\}_{n\in\N}$
defined by
$$
X_{n, t}^x = t\Phi_n(t^{-1}X_{n-1, t}^x) = t\Phi(A_n t^{-1}X_{n-1, t}^x, B_n(t^{-1}X_{n-1,
t}^x)).
$$ Then $X_{n,t}^x = t X_{n}^{t^{-1}x}$
and $\lim_{t\to 0} X_{n,t}^x = W_n^x$. Moreover, if $X_n^x$ is $\gamma$-geometric then so is $X_{n,t}^x$.
We can express $T_{t,v}$ in a slightly different form
$$
T_{t,v} f(x) = \E\Big[    e^{i\is v{X_{1,t}^x}}  f(X_{1,t}^x)\Big].
$$ 
For $t=0$ we write
$$    T_{0, v}f(x)=T_vf(x)=\E\left[e^{i\is v{W_1^x}}f(W_1^x)\right].
$$
It is not difficult to see that $h_v(x)=\E\left[e^{i\langle v, W(x)\rangle}\right]$ is an eigenfunction of $T_{v}$.
If $f\in\ct$ is an eigenfunction of operator $T_{t, v}$ with eigenvalue $k_v(t)$, then $\Dt f$ is
an eigenfunction of the operator $P_{t, v}$ with the same eigenvalue. Moreover,

\begin{lem}\label{EigenvalueP}
The unique eigenvalue of modulus $1$ for operator $P$ acting on $\ct$ is $1$ and the eigenspace is
one dimensional. The corresponding projection on $\C\cdot1$ is given by the map $f\mapsto\nu(f)$.
The unique eigenvalue of modulus $1$ for operator $T_v$ acting on $\ct$, where $v\in\Sd$, is $1$
and the eigenspace is one dimensional. The corresponding projection on $\C\cdot h_v(x)$, is given
by the map $f\mapsto f(0)\cdot h_v(x)$.
\end{lem}
\begin{proof}
For the proof, of the first part see section 3 of \cite{BDG}, and of the second part see section 5
of \cite{M1}.
\end{proof}
The lemma above says that 1 is the unique peripheral  eigenvalue both for $P_v$ and $T_v$. Even
more can be proved, the complementary part of the spectrum of both operators on $\bl1$ is contained
in a ball centered at zero and with the radius strictly smaller than 1. So, they are quasi-compact.
Moreover, due to the perturbation theorem of Keller and Liverani \cite{KL}, (see also \cite{Liv})
for small values of $t$, spectral properties of $P_{t,v}$ ($T_{t,v}$, resp.) approximate
appropriate properties of $P_{v}$ ($T_v$ resp.). The proof is based on
$\gamma$-geometricity of Markov processes $X_n^x$ and $X_{n,t}^x$, and the boundedness of $B^2$
(see Theorem \ref{thm: LT}) which in turn allows us to show that
\begin{align}\label{estim1}
|\Phi(A x, tB(t^{-1}x))-\overline{\Phi}(x)|\le \mbox{Lip}_{\Phi}|tB(t^{-1}x)|\le t \mbox{Lip}_{\Phi}(|B^1|+C),
\end{align}
for every $x\in\R^d$ and $t>0$, where $\mbox{Lip}_{\Phi}$ is the Lipschitz constant of $\Phi$.

We will not present the details, since the proof is a straightforward application of the arguments presented in \cite{BDG,M1}.

The following proposition summarizes the necessary spectral properties of operators $P_{t, v}$ and $T_{t,
v}$.
\begin{prop}\label{KL}
Assume that $0<\epsilon<1$, $\l>0$, $\l+2\epsilon<\rho=2\l$ and $2\l+\epsilon<\a$, then there exist $\d>0$,
$0<\varrho<1-\d$ and $t_0>0$ such that for every $t\in[0, t_0]$ and every
$v\in\Sd$
\begin{itemize}
\item $\sigma(P_{t, v})$ and $\sigma(T_{t, v})$ are contained in $\mathcal{D}=\{z\in\C: |z|\leq\varrho\}\cup\{z\in\C: |z-1|\leq\d\}$.
\item The sets $\sigma(P_{t, v})\cap\{z\in\C: |z-1|\leq\d\}$ and $\sigma(T_{t, v})\cap\{z\in\C:
|z-1|\leq\d\}$ consist of exactly one eigenvalue $k_v(t)$, where $\lim_{t\rightarrow0}k_v(t)=1,$
and the corresponding eigenspace is one dimensional.
\item We can express operators $P_{t,v}$ and $T_{t,v}$ in the
following form
\begin{align*}
P_{t,v}^n=k_v(t)^n\Pi_{P, t}+Q_{P, t}^n,\ \mbox{and}\ \ T_{t,v}^n=k_v(t)^n\Pi_{T, t}+Q_{T, t}^n,
\end{align*}
for every $n\in\N$,  $\Pi_{P, t}$ and $\Pi_{T, t}$ being the projections onto the one dimensional eigenspaces mentioned above. $Q_{P, t}$ and $Q_{T, t}$ are complementary operators to projections
$\Pi_{P, t}$ and $\Pi_{T, t}$ respectively, such that $\Pi_{P, t}Q_{P, t}=Q_{P, t}\Pi_{P, t}=0$ and
$\Pi_{T, t}Q_{T, t}=Q_{T, t}\Pi_{T, t}=0$. Furthermore $\|Q_{P, t}^n\|_{\bl1}=O(\varrho^n)$ and
$\|Q_{T, t}^n\|_{\bl1}=O(\varrho^n)$ for every $n\in\N$. The operators $\Pi_{P, t}$, $\Pi_{T, t}$,
$Q_{P, t}$ and $Q_{T, t}$ depend on $v\in\Sd$, but this is omitted for simplicity.
\end{itemize}
\end{prop}

The following theorem contains the basic estimate
\begin{thm}\label{diffPI}
Let $h_v$ be the eigenfunction for operator $T_v$, and assumptions of Proposition \ref{KL} are
satisfied. Then for any $0<\d\leq1$ such that $\epsilon<\d<\a$, there exist $C>0$ such that for
every $0<t\leq t_0$ we have
\begin{align}
\label{diffPI1} \left\|\D_t(\Pi_{T, t}-\Pi_{T, 0})h_v\right\|_{\rel}&\leq Ct^{\d},\ \ \ \mbox{and}\\
\label{rate1} \nu(\D_t\Pi_{T, t}h_v-1)&\leq Dt^{\d}.
\end{align}
\end{thm}
\begin{proof}
The estimate \eqref{diffPI1} bases on the inequality \eqref{estim1} and spectral properties of the
operators $T_{t,v}$. For more details we refer to section 6 in \cite{M1}.
\end{proof}
The following lemma was proved in \cite{M1} as a straightforward consequence of inequality
\eqref{diffPI1}:
\begin{lem}\label{dyaprop}
If $\a\in(0, 2)$, assumptions of Proposition \ref{KL} are satisfied and $\a-\rho> 1$ if $\a>1$, then
\begin{align}\label{dyaprop1}
\lim_{t\rightarrow0}\frac{L_b(t^{-1})}{t^{\a}}\int_{\R^d}\left(e^{it\langle v,
x\rangle}-1\right)\left(\Pi_{T, t}(h_v)(tx)-\Pi_{T, 0}(h_v)(tx)\right)\nu(dx)=0.
\end{align}
\end{lem}


\begin{proof}[Proof of Theorem \ref{thm: LT}]
Notice that $\Dt\Pi_{T, t}(h_v)$ is an eigenfunction of the operator $P_{t, v}$ corresponding to
the eigenvalue $k_v(t)$ and we have
\begin{align}
\label{Taylorf}(k_v(t)-1)\cdot\nu(\Dt\Pi_{T, t}h_v)=\nu\left(\big(e^{it\is v\cdot
}-1\big)\cdot(\Dt\Pi_{T, t}h_v)\right).
\end{align}
We will often use Theorem \ref{thm: ogon}, but in a stronger version. Observe that the limit
\begin{align}\label{con}
    \lim_{t\to0}\frac{L_b(t^{-1})}{t^{\a}}\int_{\R^d}f(tx)\nu(dx)=\L^1(f),
\end{align}
exists for every $f\in\F$, where
\begin{equation}
\label{F}
\F=\big\{f: \sup_{x\in\R^d}|x|^{-\a}|\log|x||^{1+\eps}|f(x)|<\8\ \mbox{for some $\eps>0$ and }
 \L^1(\Dis(f))=0
\big\}.
\end{equation}
Now we consider each case separately.

\bigskip

\noindent \underline{\textsc{Case $0<\a<1$.}} Observe that $\lim_{t\to0}\nu(\Dt\Pi_{T, t}h_v)=1$ by \eqref{rate1}, hence using \eqref{Taylorf} we will prove
\begin{align}\label{caseto1}
\lim_{t\rightarrow0}L_b(t^{-1})\frac{k_v(t)-1}{t^{\a}}=\int_{\R^d}\left(e^{i\langle v,
x\rangle}-1\right)h_v(x)\L^1(dx)=:C_{\a}(v).
\end{align}
Let us write,
\begin{align*}
\frac{L_b(t^{-1})}{t^{\a}}\int_{\R^d}\left(e^{it\langle v, x\rangle}-1\right)&\Pi_{T, t}(h_v)(tx)\nu(dx)\\
&=\frac{L_b(t^{-1})}{t^{\a}}\int_{\R^d}\left(e^{it\langle v, x\rangle}-1\right)\cdot\left(\Pi_{T,
t}(h_v)(tx)-\Pi_{T, 0}(h_v)(tx)\right)\nu(dx)\\
&+\frac{L_b(t^{-1})}{t^{\a}}\int_{\R^d}\left(e^{it\langle v, x\rangle}-1\right)\Pi_{T,
0}(h_v)(tx)\nu(dx).
\end{align*}
In view of Lemma \ref{dyaprop} the first term of the sum above tends to
$0$. Observe that the function $f_v(x)=\left(e^{i\langle v,
x\rangle}-1\right)h_v(x)$ belongs to $\F$ since it  is bounded
and  $|f_v(x)|\leq2|x|$ for $|x|<1$. Therefore, by Lemma \ref{dyaprop}
and (\ref{con}) the  expression above tends to a constant
as $t$ goes to 0. Thus in view of \eqref{Taylorf} we obtain \eqref{caseto1}.
Now we will show
\begin{align}\label{chcaseto1}
\lim_{n\rightarrow\8}\Xi_{\a}^n(tv)=\ups_{\a}(tv),
\end{align} where $\Xi_\a^n$ is the characteristic function of $a_n^{-1}S_n^x-d_n$.
For $t_n = \frac t{a_n}$ notice that
\begin{align*}
\Xi_{\a}^n(tv)=\E\left(e^{it_n\langle v, S_n^x\rangle}\right)=\left(P^n_{t_n,
v}1\right)(x)=k^n_v(t_n)\left(\Pi_{P, t_n}1\right)(x)+\left(Q^n_{P, t_n}1\right)(x).
\end{align*}
Since $\lim_{n\rightarrow\8}\|Q^n_{P, t_n}\|_{\bl1}=0$, by Proposition \ref{KL} and $\lim_{n\rightarrow\8}\Pi_{P, t_n}1=1$, (see \cite{BDG} or \cite{M1} for more details), we have
\begin{align*}
\lim_{n\to\8} \Xi_\a^n(tv) = \lim_{n\rightarrow\8}k^n_v(t_n)= e^{\lim_{n\to\8} n(k_v(t_n)-1)},
\end{align*}
and finally by \eqref{an} and \eqref{caseto1}
$$\lim_{n\rightarrow\8}n\cdot\left(k_v(t_n)-1\right)=\lim_{n\rightarrow\8}\frac{n\cdot t_n^{\a}}{L_b(t_n^{-1})}L_b(t_n^{-1})\frac{k_v(t_n)-1}{t_n^{\a}}
=\frac{t^{\a}C_{\a}(v)}{c}.
$$
This proves the pointwise convergence $\Xi_{\a}^n$ to $\ups_{\a}$. Continuity of $\ups_{\a}$ at $0$
follows from the Lebesgue dominated convergence theorem.


\bigskip

\noindent \underline{\textsc{Case $\a=1$.}} We prove the following
lemma
\begin{lem}\label{case1A}
For every $0<\d<1$, there exists a constant $C_{\d}>0$  such that for every $|t|\le1$,
$$|\xi(t)|
\leq C_{\d}|t|^{\d}.$$
\end{lem}
\begin{proof}
For $|t|\le 1$, we write
\begin{align*}
|\xi(t)|\leq
\int_{\R^d}\frac{|tx|}{1+|tx|^2}\nu(dx)=\left(\int_{A_1}+\int_{A_2}+\int_{A_3}\right)\left(\frac{|tx|}{1+|tx|^2}\right)\nu(dx),
\end{align*}
where  $A_1=\{x\in\R^d: |x|\leq1\}$, $A_2=\{x\in\R^d: 1<|x|\leq\frac{1}{|t|}\}$ and
$A_3=\{x\in\R^d: |x|>\frac{1}{|t|}\}$. The first integral is dominated by $C|t|$.
To estimate the third one notice that
since $\frac{|x|}{1+|x|^2}\mathbf{1}_{\{|x|>1\}}\in\F$, we have
\begin{align*}
\lim_{t\to 0}\frac{L_b(t^{-1})}{|t|}\int_{\R^d}\frac{|tx|}{1+|tx|^2}\mathbf{1}_{\{
|x|\geq\frac{1}{t}\}}\nu(dx)= \int_{\{|x|>1\}}\frac{|x|}{1+|x|^2}\L^1(dx).
\end{align*}
Therefore
$$\int_{A_3}\frac{|tx|}{1+|tx|^2}\nu(dx)\le \frac{C|t|}{L_b(t^{-1})}
\leq C|t|^{\d}.$$
Finally we estimate the second integral. Let
 $\d<\d_1<1$ and notice that $\frac 1{L_b}$ is also a slowly varying function. Then
\begin{multline*}
\sum_{k=0}^{|\log_2 |t||}\int_{\R^d}\frac{|tx|}{1+|tx|^2}\mathbf{1}_{\{
2^k<|x|\leq 2^{k+1}\}}\nu(dx)
\leq |t|\sum_{k=0}^{|\log_2 |t||}2^{k+1}\nu(\{x\in\R^d: |x|>2^k\})\\
\leq C|t|\sum_{k=0}^{|\log_2 |t||}\frac{1}{L_b(2^k)}
\le C|t|\sum_{k=0}^{|\log_2 |t||}2^{(1-\d_1)k}\le
C|t|^{\d},
\end{multline*}
since $\frac{1}{L_b(2^k)}\le C2^{(1-\d_1)k}$ (see \cite{BGT} Proposition 1.3.6 (v)).
This completes the proof of the lemma.
\end{proof}
In order to prove
\begin{align}
\label{case1}\lim_{t\rightarrow0}L_b(t^{-1})\frac{k_v(t)-1-i\langle v,
\xi(t)\rangle}{t}=\int_{\R^d}\left(\left(e^{i\langle v,
x\rangle}-1\right)h_v(x)-\frac{i\langle v, x\rangle}{1+|x|^2}\right)\L^1(dx)=:\widetilde{C}_{1}(v),
\end{align}
notice that,
\begin{align*}
\int_{\R^d}\left(e^{it\langle v, x\rangle}-1\right)\Pi_{T,
t}(h_v)(tx)\nu(dx)&=\int_{\R^d}\left(e^{it\langle v, x\rangle}-1\right)\cdot\left(\Pi_{T,
t}(h_v)(tx)-\Pi_{T, 0}(h_v)(tx)\right)\nu(dx)\\
&+\int_{\R^d}\left(e^{it\langle v, x\rangle}-1\right)\cdot\left(\Pi_{T,
0}(h_v)(tx)-1\right)\nu(dx)\\
&+\int_{\R^d}\left(e^{it\langle v, x\rangle}-1-\frac{i\langle v,
tx\rangle}{1+|tx|^2}\right)\nu(dx)+i\langle v, \xi(t)\rangle .
\end{align*}
The first term of the sum tends to $0$ by Lemma \ref{dyaprop}. The
function $f_v(x)=\left(e^{i\langle v,
x\rangle}-1\right)(h_v(x)-1)$ belongs to $\F$. Indeed, $f_v$ is
bounded and for $|x|<1$
\begin{align*}
|f_v(x)|\leq\left|e^{i\langle x, v\rangle}-1\right||h_v(x)-1|\leq2\E(|W(x)|^{\d})|x|\le C|x|^{1+\d},
\end{align*}
for any $0<\d<1$. Similarly, one can prove that $g_v(x)=e^{i\langle v, x\rangle}-1-\frac{i\langle v,
x\rangle}{1+|x|^2}$ belongs to $\F$. Indeed, $g_v$ is bounded and for $|x|<1$
\begin{align*}
|g_v(x)|\leq\left|e^{i\langle x, v\rangle}-1-i\langle x,
v\rangle\right|+\frac{|x|^3}{1+|x|^2}\leq2|x|^{1+\d}+\frac{|x|^3}{1+|x|^2},
\end{align*}
for any $0<\d<1$. Hence, by (\ref{con}) we obtain
\begin{align}
\label{case1b}\lim_{t\rightarrow0}\frac{L_b(t^{-1})}{t}\left(\int_{\R^d}\left(e^{it\langle v,
x\rangle}-1\right)\Pi_{T, t}(h_v)(tx)\nu(dx)-i\langle v, \xi(t)\rangle\right)=\widetilde{C}_{1}(v).
\end{align}
Now by (\ref{case1b}) we have
\begin{multline*}
\lim_{t\rightarrow0}L_b(t^{-1})\frac{k_v(t)\!-\!1\!-\!i\langle v,
\xi(t)\rangle}{t}
=\lim_{t\rightarrow0}L_b(t^{-1})\frac{\left(\nu\left((e^{it\langle v,
\cdot\rangle}\!\!-\!1)(\Dt\Pi_{T, t}h_v)\right)\!\!-\!i\langle v, \xi(t)\rangle\nu(\Dt\Pi_{T,
t}h_v)\right)}{\nu(\Dt\Pi_{T, t}h_v)t}\\
=\lim_{t\rightarrow0}L_b(t^{-1})\left(\frac{\nu\left((e^{it\langle v,
\cdot\rangle}-1)(\Dt\Pi_{T, t}h_v)\right)-i\langle v, \xi(t)\rangle}{\nu(\Dt\Pi_{T,
t}h_v)t}+\frac{i\left(1-\nu(\Dt\Pi_{T, t}h_v)\right)\langle v, \xi(t)\rangle}{\nu(\Dt\Pi_{T,
t}h_v)t}\right)
=\widetilde{C}_1(v).\\
\end{multline*}
Since by \eqref{rate1} and Lemma \ref{case1A} we have
\begin{align*}
\lim_{t\rightarrow0}L_b(t^{-1})\left(\frac{i\left(1-\nu(\Dt\Pi_{T, t}h_v)\right)\langle v,
\xi(t)\rangle}{\nu(\Dt\Pi_{T, t}h_v)t}\right)=0,
\end{align*}
and (\ref{case1}) follows. Now we need the following

\begin{lem}
Let $m_{\sigma_{\L^1}}=\int_{\mathbb{S}^{d-1}}\o\sigma_{\L^1}(d\o)$, where $\sigma_{\L^1}$ is the spherical measure associated with the tail measure $\L^1$. Then for every $t\in\R$ and $v\in\Sd$
\begin{align}
\label{chcase1a}\lim_{s\rightarrow0}\frac{L_b(s^{-1})}{s}\int_{\R^d}\left(\frac{\langle v,
stx\rangle}{1+|stx|^2}-\frac{\langle v, stx\rangle}{1+|sx|^2}\right)\nu(dx)=-t\log|t|\langle v, m_{\sigma_{\L^1}}\rangle.
\end{align}
In particular, for every $0<\d<1$ there exists a constant $C_{\d}>0$  such that for every $|t|\le1$,
\begin{align}
\label{chcase2aa}|t\log|t|\langle v, m_{\sigma_{\L^1}}\rangle|\leq C_{\d} |t|^{\d}.
\end{align}
\end{lem}
\begin{proof}
Observe that $\frac{x}{1+|tx|^2}-\frac{x}{1+|x|^2}\in\mathcal{F}$ hence
 \begin{align*}
\lim_{s\rightarrow0}\frac{L_b(s^{-1})}{s}\int_{\R^d}\left(\frac{\langle v,
stx\rangle}{1+|stx|^2}-\frac{\langle v, stx\rangle}{1+|sx|^2}\right)\nu(dx)=t\langle v,
\tau(t)\rangle,
\end{align*}
where $\tau(t)=\int_{\R^d}\left(\frac{x}{1+|tx|^2}-\frac{x}{1+|x|^2}\right)\L^1(dx)$.
Notice that
\begin{multline*}
    \tau(t)=\int_{\R^d}\left(\frac{x}{1+|tx|^2}-\frac{x}{1+|x|^2}\right)\L^1(dx)
    =\int_0^{\8}\int_{\mathbb{S}^{d-1}}
    \left(\frac{r\o}{1+|tr\o|^2}-\frac{r\o}{1+|r\o|^2}\right)\sigma_{\L^1}(d\o)\frac{dr}{r^2}\\
    =\int_{\mathbb{S}^{d-1}}\o\sigma_{\L^1}(d\o)\cdot\int_0^{\8}
    \left(\frac{r(1-t^2)}{(1+t^2r^2)(1+r^2)}\right)dr
    =-m_{\sigma_{\L^1}}\log|t|.
\end{multline*}
The proof is completed.
\end{proof}


%
For $t_n = \frac t{a_n}$, $t>0$ we have
\begin{multline*}
\lim _{n\to\8}\Xi_{1}^n(tv) 
=\lim_{n\to\8} e^{-itn\langle v,\xi(a_n^{-1})\rangle}\E\left(e^{it_n\langle v, S_n^x\rangle}\right)\\
=\lim_{n\rightarrow\8}e^{-int\langle v, \xi(a_n^{-1})\rangle}k_v^n(t_n)= e^{\lim_{n\to\8}
\big(n(e^{-it\langle v, \xi(a_n^{-1})\rangle}k_v(t_n)-1)\big)}.
\end{multline*}
Hence
\begin{multline*}
\lim_{n\rightarrow\8}\left(n\left(e^{-it\langle v,
\xi(a_n^{-1})\rangle}k_v(t_n)-1\right)\right)\\
=\lim_{n\rightarrow\8}\!\!\left(\frac{nt_n}{L_b(t_n^{-1})}e^{-it\langle v, \xi(a_n^{-1})\rangle}L_b(t_n^{-1})\frac{k_v(t_n)-1-i\langle v, \xi(t_n)\rangle}{t_n}+ ne^{-it\langle v, \xi(a_n^{-1})\rangle}\left(1+i\langle v,
\xi(t_n)\rangle\right)-n\right)\\
=\lim_{n\rightarrow\8}\!\!
\left(\widetilde{C}_{1}(v)\frac{nt}{a_nL_b(a_n)}\frac{L_b(a_n)}{L_b(t^{-1}a_n)}\!+\!
 n\left(1\!-\!it\langle v, \xi(a_n^{-1})\rangle+O\left(t^2\langle v, \xi(a_n^{-1})\rangle^2\right)\right)\left(1+i\langle v, \xi(t_n)\rangle\right)-n\right)\\
=\lim_{n\rightarrow\8}\!\!\Big(in\langle v, \xi(t_n)\rangle-int\langle
v,
\xi(a_n^{-1})\rangle+
nt\langle v, \xi(t_n)\rangle\langle v,
\xi(a_n^{-1})\rangle+nO\left(t^2\langle v, \xi(a_n^{-1})\rangle^2\right)\left(1+i\langle v,
\xi(t_n)\rangle\right)\Big)\\
+ \frac{t\widetilde{C}_{1}(v)}{c}.
\end{multline*}
Notice that by (\ref{chcase1a}) we have
\begin{multline*}
\lim_{n\rightarrow\8}\left(in\langle v, \xi(t_n)\rangle-int\langle v,
\xi(a_n^{-1})\rangle\right)=\\
=\lim_{n\rightarrow\8}it\frac{n}{a_nL_b(a_n)}\cdot a_nL_b(a_n)\int_{\R^d} \left(\frac{\langle v,
a_n^{-1}x\rangle}{1+\left|a_n^{-1}tx\right|^2}-\frac{\langle v,
a_n^{-1}x\rangle}{1+\left|a_n^{-1}x\right|^2}\right)\nu(dx)\\
=-\frac{it\log t\langle v, m_{\sigma_{\L^1}}\rangle}{c},
\end{multline*}
and the limit of two remaining factors, by Lemma
\ref{case1A}, is 0. Therefore the limit of the whole expression
 is equal to
$\ups_1(tv)=\frac{t\widetilde{C}_1(v)-it\log t\langle v, m_{\sigma_{\L^1}}\rangle}{c}$. Finally, to prove continuity of $\ups_1$
at zero, it is enough to observe that for $|x|<1$,
\begin{align*}
\left|\left(e^{i\langle v, x\rangle}-1\right)h_v(x)-\frac{i\langle v, x\rangle}{1+|x|^2}\right|\leq
C|x|^{1+\d},
\end{align*}
for any $0<\d<1$ and some $C>0$ independent of $v\in\Sd$.


\bigskip

\noindent \underline{\textsc{Case $1<\a<2$.}} As in the previous cases we show that
\begin{align}
\label{case1to2}\lim_{t\rightarrow0}L_b(t^{-1})\frac{k_v(t)-1-i\langle v,
tm\rangle}{t^{\a}}=\int_{\R^d}\left(\left(e^{i\langle v, x\rangle}-1\right)h_v(x)-i\langle v,
x\rangle\right)\L^1(dx)=:C_{\a}(v).
\end{align}

Let us write,
\begin{align*}
\int_{\R^d}\left(e^{it\langle v, x\rangle}-1\right)\Pi_{T,
t}(h_v)(tx)\nu(dx)&=\int_{\R^d}\left(e^{it\langle v, x\rangle}-1\right)\cdot\left(\Pi_{T,
t}(h_v)(tx)-\Pi_{T, 0}(h_v)(tx)\right)\nu(dx)\\
&+\int_{\R^d}\left(e^{it\langle v, x\rangle}-1\right)\cdot\left(\Pi_{T,
0}(h_v)(tx)-1\right)\nu(dx)\\
&+\int_{\R^d}\left(e^{it\langle v, x\rangle}-1-i\langle v, tm\rangle\right)\nu(dx)+i\langle v,
tm\rangle,
\end{align*}
By Lemma \ref{dyaprop} the first term of the sum goes to $0$. Functions
$f_v(x)=\left(e^{i\langle v, x\rangle}-1\right)(h_v(x)-1)$ and
$g_v(x)=e^{i\langle v, x\rangle}-1-i\langle v, x\rangle$  belong
to $\F$.
Hence, by (\ref{con}) we obtain
\begin{align}
\label{case1to2a}\lim_{t\rightarrow0}\frac{L_b(t^{-1})}{t^{\a}}\left(\int_{\R^d}\left(e^{it\langle v,
x\rangle}-1\right)\Pi_{T, t}(h_v)(tx)\nu(dx)-i\langle v, tm\rangle\right)=C_{\a}(v).
\end{align}
Similarly, as in the previous case we have
\begin{multline*}
\lim_{t\rightarrow0}L_b(t^{-1})\frac{k_v(t)-1-i\langle v,
tm\rangle}{t^{\a}}=\lim_{t\rightarrow0}L_b(t^{-1})\frac{\left(\nu\left((e^{it\langle v,
\cdot\rangle}\!-\!1)(\Dt\Pi_{T, t}h_v)\right)\!-\! i\langle v, tm\rangle\nu(\Dt\Pi_{T,
t}h_v)\right)}{\nu(\Dt\Pi_{T, t}h_v)t^{\a}}\\
\lim_{t\rightarrow0}L_b(t^{-1})\left(\frac{\nu\left((e^{it\langle v, \cdot\rangle}-1)\cdot(\Dt\Pi_{T,
t}h_v)\right)-i\langle v, tm\rangle}{\nu(\Dt\Pi_{T, t}h_v)t^{\a}}+\frac{i\left(1-\nu(\Dt\Pi_{T,
t}h_v)\right)\langle v, tm\rangle}{\nu(\Dt\Pi_{T, t}h_v)t^{\a}}\right)=C_{\a}(v).
\end{multline*}
Since by \eqref{rate1}
\begin{align*}
\lim_{t\rightarrow0}L_b(t^{-1})\left(\frac{i\left(1-\nu(\Dt\Pi_{T, t}h_v)\right)\langle v,
tm\rangle}{\nu(\Dt\Pi_{T, t}h_v)t^{\a}}\right)=0,
\end{align*}
and (\ref{case1to2}) follows.

Now we can show that
\begin{align}\label{chcase1to2}
\lim_{n\rightarrow\8}\Xi_{\a}^n(tv)=\ups_{\a}(tv),
\end{align}
In order to prove (\ref{chcase1to2}) notice that
$$
\lim_{n\to\8}\Xi_{\a}^n(tv)= \lim_{n\to\8}e^{-int_n\langle v, m\rangle}\E\left(e^{it_n\langle v,
S_n^x\rangle}\right)
=e^{\lim_{n\rightarrow\8} n\left(e^{-it_n\langle v,
m\rangle}k_v(t_n)-1\right)}.
$$
Moreover, since $\lim_{n\rightarrow\8}nt_n^2 = 0$, we have
\begin{multline*}
\lim_{n\rightarrow\8}\left(n\left(e^{-it_n\langle v,
m\rangle}k_v(t_n)-1\right)\right)\\
=\lim_{n\rightarrow\8}\left(\frac{nt_n^{\a}}{L_b(t_n^{-1})}e^{-it_n\langle v,
m\rangle}\cdot L_b(t_n^{-1})\frac{k_v(t_n)-1-it_n\langle v, m\rangle}{t_n^{\a}}+ne^{-it_n\langle v,
m\rangle}\left(1+it_n\langle v,
m\rangle\right)-n\right)\\
=\lim_{n\rightarrow\8} \left( C_{\a}(v)\frac{nt^{\a}}{a_n^{\a}L_b(a_n)} \frac{L_b(a_n)}{L_b(t^{-1}a_n)}
+\left(n\cdot\left(1-it_n\langle v,
m\rangle+O\left(t_n^2\right)\right)\cdot\left(1+it_n\langle v, m\rangle\right)-n\right)\right)\\
=\frac{t^{\a}C_{\a}(v)}{c}+\lim_{n\rightarrow\8}\left(nt_n^2\langle v,
m\rangle^2+nO\left(t_n^2\right)\cdot\left(1+it_n\langle v,
m\rangle\right)\right)=\frac{t^{\a}C_{\a}(v)}{c},\\
\end{multline*}
and (\ref{chcase1to2})
follows. To prove continuity of $\ups_{\a}$ at zero, we proceed as in the previous cases.

\bigskip

Finally, under some additional assumptions, we have to prove a nondegeneracy of the limit variable $C_{\a}(v)$ for $v\in\Sd$.
Notice first that since $\Phi(x,0)=x$ for every $x\in\overline{[\supp\mu]\cdot\supp\nu}$,
$W(x)=\sum_{k=1}^{\8}A_{k}\cdot\ldots\cdot A_{1}x$.
Let us define $W^*(x)=\sum_{k=1}^{\8}A_{1}^*\cdot\ldots\cdot A_{k}^*x$ and observe
\begin{eqnarray*}
\Re C_\a(v)&=& \Re \bigg(
\int_{\R^d}\big( e^{i \is vx} - 1 \big)\E\big[ e^{i\is v{W(x)}} \big] \L^1(dx)
\bigg)\\
&=& \int_0^\8
\int_{\S^{d-1}}\E\Big[ \cos\big( t\is {W^*(v)+v}w \big)
- \cos\big( t\is {W^*(v)}w \big)
\Big] \sigma_{\L^1}(dw)\frac {dt}{t^{\a+1}}.
\end{eqnarray*}
Hence
$$
\Re C_\a(v) = C(\a)\cdot \int_{\S^{d-1}} \E\Big[ \big| \is {W^*(v)+v}w \big|^\a - \big| \is {W^*(v)}w
\big|^\a \Big]\sigma_{\Lambda^1}(dw),
$$ for $C(\a) = \int_0^\8 \frac{\cos t-1}{t^{\a+1}}dt<0$. Notice that $W_v = W^*(v)+v$ is a solution
of the random difference equation
\begin{equation}
\label{rde}
W_v =_d A^* W_v + v.
\end{equation}
Moreover, since $\lim_{n\to\8}\left(\E\|A_1\cdot\ldots\cdot A_n\|^{\a}\right)^{\frac 1 n}<1$, implies that $\E|W_v|^\a<\8$, we have
\begin{align*}
\Re C_\a(v) &= C(\a)\cdot \int_{\S^{d-1}} \E\Big[ \big| \is {W_v}w \big|^\a - \big| \is {A^*W_v}w
\big|^\a \Big]\sigma_{\Lambda^1}(dw),\\
&=C(\a)\cdot \int_{\S^{d-1}} \E\Big[ \big| \is {W_v}w \big|^\a - |Aw|^{\a}\big| \is {W_v}{A* w}
\big|^\a \Big]\sigma_{\Lambda^1}(dw).
\end{align*}
Now in view of \eqref{ogon1b} and \eqref{ogon1c} we obtain
\begin{align*}
\int_{\S^{d-1}} \E\Big[|Aw|^{\a}\big| \is {W_v}{A* w}
\big|^\a \Big]\sigma_{\Lambda^1}(dw)=\sum_{n=2}^{\8}\int_{\S^{d-1}} \E\Big[\is {W_v}{w}
\big|^\a \Big]\sigma_{\Gamma_n}(dw).
\end{align*}
Therefore 
we can conclude that for every $v\in\mathbb{S}^{d-1}$
$$\Re C_\a(v) = C(\a)\cdot \int_{\S^{d-1}} \E\Big[ \big| \is {W_v}w \big|^\a \Big]\sigma_{\Gamma_1}(dw).$$
Finally we have to prove that $\int_{\S^{d-1}} \E\Big[ \big| \is {W_v}w \big|^\a \Big]\sigma_{\Gamma_1}(dw)>0$. For this purpose, in view of \eqref{eq:lambda1}, notice that for every $f\in C(\S^{d-1})$
$$\int_{\S^{d-1}}f(w)\sigma_{\Gamma_1}(dw)=\int_{\S^{d-1}}f\left(\frac{\Phi(0, w)}{|\Phi(0, w)|}\right)|\Phi(0, w)|^{\a}\sigma_{\L_b}(dw),$$
which in turn implies
\begin{multline}\label{degener1}
    \int_{\S^{d-1}} \E\Big[ \big| \is {W_v}w \big|^\a \Big]\sigma_{\Gamma_1}(dw)=\int_{\S^{d-1}} \E\Big[ \big| \is {W_v}{\Phi(0, w)} \big|^\a \Big]\sigma_{\Lambda_b}(dw)\\
    =\E\Bigg[|W_v|^{\a}\int_{\S^{d-1}}  \big| \is {W_v/|W_v|}{\Phi(0, w)} \big|^\a \sigma_{\Lambda_b}(dw)\Bigg]\ge C_{\L_b}\E\big[|W_v|^{\a}\big]
\end{multline}
for $C_{\L_b}=\min_{u\in\S^{d-1}}\int_{\S^{d-1}}  \big| \is {u}{\Phi(0, w)} \big|^\a \sigma_{\Lambda_b}(dw)$  which is strictly positive. Indeed, if for some $u_0\in\S^{d-1}$,  $\int_{\S^{d-1}}  \big| \is {u_0}{\Phi(0, w)} \big|^\a \sigma_{\Lambda_b}(dw)$ were equal to 0, then
 the set $\Phi[\{0\}\times\supp\sigma_{\Lambda_b}]$ would be contained in the hyperplane $u_0^{\bot}$, that contradicts to our assumptions.
This completes the proof of Theorem \ref{thm: LT}.
\end{proof}


\bibliographystyle{alpha}
\newcommand{\etalchar}[1]{$^{#1}$}

\end{document}